\documentclass[final]{siamltex}

\usepackage{latexsym, amsmath,  amssymb}
\usepackage{graphics}
\usepackage{colortbl}
\usepackage{graphicx}

\newcommand{\A}{\mathbf A}
\newcommand{\B}{\mathbf B}
\newcommand{\HH}{\mathbb H^2}
\newcommand{\bb}{\mathbf b}
\newcommand{\nn}{\mathbf n}
\newcommand{\pp}{\mathbf e_1}
\newcommand{\qq}{\mathbf e_2}

\newcommand{\T}{\mathbf T}
\newcommand{\Ts}{\T_s}
\newcommand{\Tss}{\T_{ss}}
\newcommand{\Xt}{\mathbf X_t}
\newcommand{\X}{\mathbf X}
\newcommand{\Xs}{\mathbf X_s}
\newcommand{\Xss}{\mathbf X_{ss}}
\newcommand{\Tt}{\mathbf T_t}

\newcommand{\propn}{\wedge_\pm}
\newcommand{\prop}{\wedge_{+}} 
\newcommand{\pron}{\wedge_{-}} 

\newcommand{\cirpn}{\circ_\pm}

\title{\textbf{A numerical study of the self-similar solutions of the
Schr\"odinger Map}}

\author{F. de la Hoz\thanks{Departamento de Matem\'atica Aplicada,
    Universidad del Pa{\'\i}s Vasco-Euskal Herriko Unibertsitatea.
    Partially supported by grant MTM2007-62186} \and C.
    Garc{\'\i}a-Cervera\thanks{Mathematics Department,
    University of California,
    Santa Barbara, CA 93106, USA. Partially supported by NSF grant DMS-0505738.},
    \and L. Vega\thanks{Departamento de Matem\'aticas,
    Universidad del Pa{\'\i}s Vasco-Euskal Herriko Unibertsitatea.
    Partially supported by grant MTM2007-62186}}

\begin{document}

\maketitle

\begin{abstract}
{\sl We present a numerical study of  the self-similar solutions of the
Localized Induction Approximation of a vortex filament. These
self-similar solutions, which constitute a one-parameter family,
develop a singularity at finite time. We study a number of boundary
conditions that allow us reproduce the  mechanism of singularity
formation. Some related questions are also considered.}
\end{abstract}



\pagestyle{myheadings} \thispagestyle{plain} \markboth{Francisco de la
Hoz, Carlos Garc{\'\i}a-Cervera and Luis Vega}{A numerical study of the
self-similar solutions of the Schr\"odinger Map}

\section{Introduction}
Given a curve $\X_0 : \mathbb{R} \longrightarrow \mathbb{R}^3$,
consider the geometric flow
\begin{equation}
\Xt = c\bb,
\end{equation}
where $c$ is the curvature and $\bb$ the binormal component of the
Frenet-Serret formulae
\begin{align}
\label{frenet-serret}
\begin{pmatrix}
\T \cr \mathbf n \cr \mathbf b
\end{pmatrix}_s =
\begin{pmatrix}
0 & c & 0 \cr -c & 0 & \tau \cr 0 & -\tau & 0
\end{pmatrix} \cdot
\begin{pmatrix}
\T \cr \mathbf n \cr \mathbf b
\end{pmatrix}.
\end{align}

\noindent The flow can be expressed as
\begin{equation}
\label{xtpos}\Xt = \Xs\prop\Xss,
\end{equation}
where $\prop$ is the usual cross-product, and $s$ denotes arclength.
This flow appeared for the first time in 1906 \cite{darios} and was
rederived in 1965 by Arms and Hama \cite{arms} as an approximation
of the dynamics of a vortex filament under the Euler equations. This
model is usually known as the Localized Induction Approximation
(LIA). We refer the reader to \cite{batchelor} and \cite{saffman}
for an analysis and discussion about the limitations of this model.
Starting with the work of Schwartz in \cite{schwarz}, the LIA has
been also used as an approximation of the quantum vortex motion in
superfluid Helium. Of particular relevance for our purposes  is the
recent work of T. Lipniacki \cite{lipniacki1}, \cite{lipniacki2}. A
rather complete list of references about the use of LIA in this
setting can be found in these two papers.\par

Some of the explicit solutions of \eqref{xtpos} are the line, circle
and helix. It is easy to see that the tangent vector $\T = \X_s$
remains with constant length, so that we can assume that it takes
values on the unit sphere. Differentiating \eqref{xtpos}, we get the
following equation for $\T$:
\begin{equation}
\label{ttpos}\T_t = \T\wedge_+\Tss.
\end{equation}
This equation, known as the Schr\"odinger map equation on the
sphere, is a particular case of the Landau-Lifshitz equation for
ferromagnetism \cite{landau} and can be rewritten in a more
geometric way as
\begin{equation}
\T_t = \mathbf J\mathbf D_s\mathbf \Ts,
\end{equation}
where $\mathbf D$ is the covariant derivative and $\mathbf J$ is the
complex structure of the sphere. Written in this way, \eqref{ttpos}
admits an immediate generalization and we can change both its
definition domain (considering for instance more variables) and its
image (considering other more complex varieties). We will insist on
the second possibility, choosing also the hyperbolic plane $\HH$ as
the target space; in that case, the equation for $\T$ is
\begin{equation}
\T_t = \T\wedge_-\Tss
\end{equation}
and, equivalently, for $\X$,
\begin{equation}
\X_t = \X_s\wedge_-\X_{ss},
\end{equation}
with $\pron$ defined as
$$
\mathbf a\pron\mathbf b=(a_2b_3-a_3b_2,a_3b_1-a_1b_3, -
(a_1b_2-a_2b_1)).
$$
In this article, we study numerically the self-similar solutions of
\begin{align}
\label{xtposneg} \Xt = \Xs\propn\Xss.
\end{align}

\noindent where $\propn$ has been defined as
$$
\mathbf a\propn\mathbf b=(a_2b_3-a_3b_2,a_3b_1-a_1b_3, \pm
(a_1b_2-a_2b_1)).
$$

\noindent Equivalently, a generalized version of the scalar product,
denoted by $\cirpn$, is given as
$$
\mathbf a\cirpn\mathbf b = a_1b_1 + a_2b_2 \pm a_3b_3.
$$

\noindent By using $\pm$, we can consider simultaneously the
Euclidean case, corresponding to $+$ and the Hyperbolic case,
corresponding to $-$. Equivalently, whenever we use $\mp$, the $-$
sign will refer to the Euclidean case and the $+$ sign to the
Hyperbolic one. Using the $\propn$ notation, the equation for $\Xs =
\T$ is now
\begin{align}
\label{ttposneg} \Tt = \T\propn\Tss.
\end{align}

\noindent If $\T\in\HH$, it is still possible to give a generalized
version of the Frenet-Serret trihedron \eqref{frenet-serret} for each
point of the curve $\X$, formed by $\T$ and two other vectors $\pp$ and
$\qq$. Indeed, a few calculations show that, for both cases, all the
possible generalizations of \eqref{frenet-serret} have the form
\begin{equation}
\label{triedroabg pn}
\begin{pmatrix}
\T \cr \pp \cr \qq\end{pmatrix}_s =
\begin{pmatrix}
0 & \alpha & \beta \cr \mp\alpha & 0 & \delta \cr \mp\beta & -\delta
& 0
\end{pmatrix} \cdot \begin{pmatrix}\T \cr \pp \cr
\qq\end{pmatrix},
\end{equation}

\noindent where $\T\cirpn\T = \pm 1$, $\pp\cirpn\pp = \qq\cirpn\qq =
1$, $\T\cirpn\pp = \T\cirpn\qq = \pp\cirpn\qq = 0$; with regard to
$\propn$, we have $\T\propn\pp = \qq$, $\pp\propn\qq = -\T$,
$\pp\propn\qq = \pm\T$. Without loss of generality, we can choose
one of the coefficients $\alpha$, $\beta$ or $\delta$ to equal zero.
If we make $\beta = 0$ and denote $\alpha \equiv c$ and $\delta
\equiv \tau$, \eqref{triedroabg pn} becomes
\begin{align}
\label{frenet-serretgen}
\begin{pmatrix}
\T \cr \pp \cr \qq
\end{pmatrix}_s =
\begin{pmatrix}
0 & c & 0 \cr \mp c & 0 & \tau \cr 0 & -\tau & 0
\end{pmatrix} \cdot
\begin{pmatrix}
\T \cr \pp \cr \qq
\end{pmatrix}.
\end{align}
In the Euclidean case, $\pp$ and $\qq$ refer to $\nn$ and $\bb$
respectively. In the hyperbolic case, we can refer to $c$ and $\tau$
as the generalized curvature and torsion. In both cases, it is
possible to recover $\X$, except for some rigid motion, if $c$ and
$\tau$ are known.

The self-similar solutions of \eqref{xtposneg} have been studied in
\cite{gutierrez} in the Euclidean case and in \cite{delahoz} in the
hyperbolic one -see also \cite{lipniacki1}, \cite{lipniacki2} for
some related work. They are such that if $\X(s, t)$ solves
\eqref{xtposneg}, so does $\lambda^{-1}\mathbf X(\lambda s,
\lambda^2t)$. Therefore, taking $\lambda = t^{-1/2}$ and defining
$\mathbf G(s) = \mathbf X(s, 1)$, the self-similar solutions will be
of the form
\begin{equation}
\mathbf X(s, t) = t^{1/2}\mathbf X(t^{-1/2}s, 1) = \sqrt t \mathbf
G(s / \sqrt t).
\end{equation}
Bearing this in mind, after some straightforward but tedious
computations \cite{gutierrez,delahoz}, the family of solutions in
which we are interested is defined by
\begin{equation}
\label{0xc0}\mathbf X_{c_0}(s,t) = \sqrt t\mathbf G(s / \sqrt t),
\end{equation}

\noindent where $c_0$ is the family parameter and $\mathbf G' =
\T(s, 1)$ is the solution of
\begin{align}
\label{Tpqt}
\begin{pmatrix}
\T \cr \pp \cr \qq
\end{pmatrix}_s =
\begin{pmatrix}
0 & c_0 & 0 \\ \mp c_0 & 0 & \frac{s}{2} \\ 0 & -\frac{s}{2} & 0
\end{pmatrix} \cdot
\begin{pmatrix}
\T \cr \pp \cr \qq
\end{pmatrix},
\end{align}

\noindent $\T\equiv \X_s$, $\pp$ and $\qq$ being the components of
the generalized Frenet-Serret formulae; the initial conditions for
\eqref{Tpqt} are
\begin{equation}
\label{Tpqt0}
\begin{cases}
\mathbf G(0) = 2c_0(0,1,0), \cr \T(0, 1) = (0,0,1), \cr \pp(0, 1) =
(1, 0, 0), \cr \qq(0, 1) = (0, 1, 0).
\end{cases}
\end{equation}

\noindent Finally, for $c_0 = 0$ we define
\begin{equation}
\label{0x0}\mathbf X_0(s,t) = s(0, 0, 1).
\end{equation}

\noindent It can be easily shown that, for an arbitrary $t > 0$, we
have
\begin{align}
\label{Tpqtt}
\begin{pmatrix}
\T \cr \pp \cr \qq
\end{pmatrix}_s =
\begin{pmatrix}
0 & \frac{c_0}{\sqrt t} & 0 \cr \mp \frac{c_0}{\sqrt t} & 0 &
\frac{s}{2t} \cr 0 & -\frac{s}{2t} & 0
\end{pmatrix} \cdot
\begin{pmatrix}
\T \cr \pp \cr \qq
\end{pmatrix}.
\end{align}

\noindent The one-parameter family of self-similar functions we have
just defined satisfies the following theorem:
\begin{theorem}
\label{teorema} Given $c_0 \ge 0$, the $\mathbf X_{c_0}$ defined by
\eqref{0xc0}, \eqref{Tpqt0} and \eqref{0x0} is a $C^\infty$ solution
of \eqref{xtposneg}, $\forall t>0$.

Moreover, there are $\mathbf A^1(c_0)$, $\mathbf A^2(c_0)$, $\mathbf
B^1(c_0)$, $\mathbf B^2(c_0)$ and a constant $C$ such that
\begin{enumerate}
    \item[(i)] $ |\mathbf X_{c_0}(s,t) - \mathbf A^1s(c_0)\chi_{[0,+\infty)}(s) -
    \mathbf A^2s(c_0)\chi_{(-\infty,0]}(s)| \le C\sqrt t. $

    \item[(ii)] We have the following asymptotics:
    \begin{align*}
    \mathbf G(s) & = \mathbf A^j(c_0)\left(s \pm 2 \dfrac{c_0^2}{s}\right)
    - 4c_0\dfrac{\pp}{s^2} + \mathcal O(1 / s^3), \qquad s\to\pm\infty; \\
    \T(s) & = \mathbf A^j(c_0) - 2c_0\dfrac{\qq}{s} + \mathcal O(1 / s^2),
\qquad s\to\pm\infty; \\
    (\pp - i\qq) & = \mathbf B^j(c_0)e^{is^2/4}e^{\pm ic_0^2\log|s|} + \mathcal O(1
/ s), \qquad s\to\pm\infty;
    \end{align*}
    where the superindex $j = 1$ when $s\to\infty$ and $j =
    2$ when $s\to-\infty$.

    \item[(iii)] $\mathbf A^j = (A_1^j, A_2^j, A_3^j)$ and $\mathbf B^j = (B_1^j,
    B_2^j, B_3^j)$ are vectors, with $\mathbf A^j\cirpn\mathbf A^j = \pm 1$ and
    \begin{align*}
        A_1^1 & = -A_1^2, \qquad A_2^1 = -A_2^2, \qquad A_3^1 = A_3^2 =
        e^{\mp\frac{c_0^2}{2} \pi}, \\
        B_1^1 & = B_1^2, \qquad B_2^1 = B_2^2, \qquad B_3^1 = -B_3^2,
        \qquad \mathbf A^j\circ_\pm\mathbf B^j = 0,
    \end{align*}
\end{enumerate}

\end{theorem}
\noindent We have combined in this theorem the results for both the
Euclidean and the Hyperbolic cases, proved in \cite{gutierrez} and
\cite{delahoz}, respectively. Naturally, $\mathbf A^j$, $\mathbf B^j$
and $C$ are different in each case.

It follows from this theorem that there exists a solution to
\eqref{xtposneg} such that
\begin{align}
\X(s, 0) = \mathbf A^1s\chi_{[0,+\infty)}(s) + \mathbf
A^2s\chi_{(-\infty,0]}(s),
\end{align}
and, correspondingly for  \eqref{ttposneg},
\begin{align}
\T(s, 0) = \A^1\chi_{[0,+\infty)}(s) + \A^2\chi_{(-\infty,0]}(s).
\end{align}
Observe that both \eqref{xtposneg} and \eqref{ttposneg} are time
reversible, because if $\X(s, t)$ and $\T(s, t)$ are their
respective solutions, so are $\X(-s, -t)$ and $-\T(-s, -t)$.
\begin{figure}[!ht]
\centering
\resizebox{12cm}{5cm}{\includegraphics{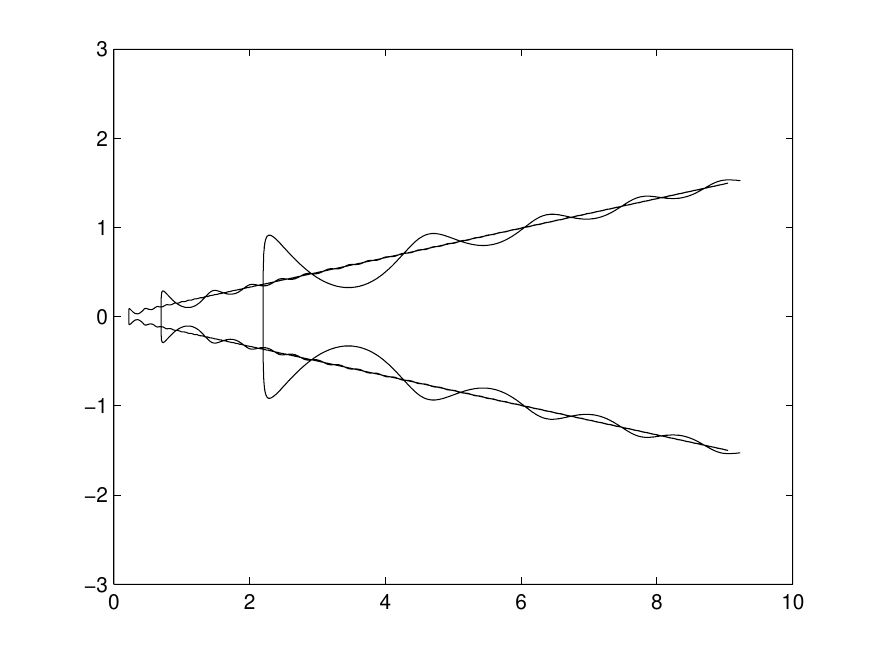}\includegraphics{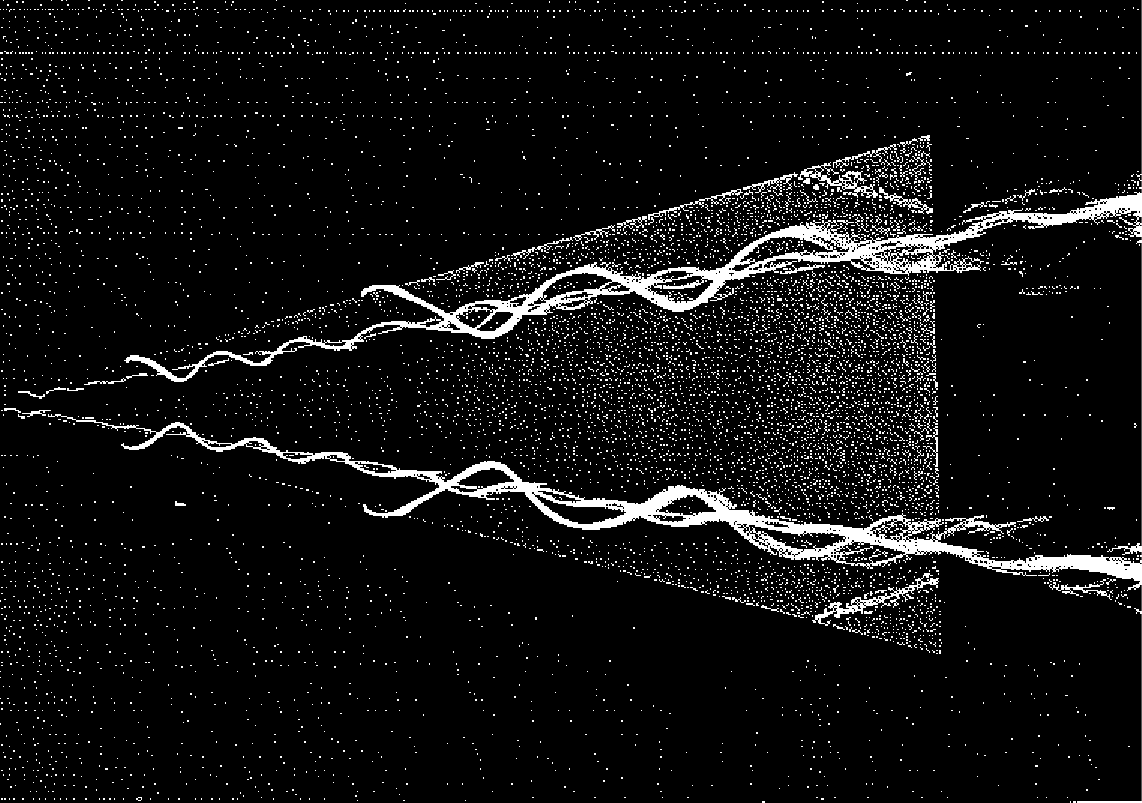}}
{\caption{\small Comparison between the theoretical evolution of
$\X(s, t)$ and a real experiment of a colored fluid traversing a
triangular wing. The experimental result on the right has been taken
from \cite{werle}.} \label{onera}}
\end{figure}
Therefore, we have two one-parameter families of regular solutions
that develop singularities at finite time. In the case of $\X$, we
have precisely a corner-shaped singularity, having thus a good model
to describe some natural phenomena. In Figure \ref{onera}, for
instance, we have plotted on the left hand-side the evolution of
$\X(s, t)$, obtained by integrating \eqref{Tpqtt}, and on the right
hand side, a picture of an experiment published in ONERA
\cite{werle}, where several lines of colored fluid in water show the
symmetrical pair of vortices behind an inclined delta wing, with a
high Reynolds number. The resemblance between both images is really
striking, at least at the qualitative level. In fact, in both images
we clearly  see the self-similar behavior of the evolution of the
filaments. Also in both cases the filaments tend asymptotically to
two non-parallel lines. Finally, both pictures have in common the
oscillatory behavior, and the shift of the horse-shoe-type of curve
that is close to the vertex with respect to the plane that contains
the two asymptotic lines. See also figures 7.5.7 and 7.8.6 of plate
2 in \cite{batchelor}. It would be interesting to know if this
resemblance is not just at the qualitative level. Notice that
Theorem \ref{teorema} quantifies in a precise manner the above
mentioned properties for the self-similar solutions of the LIA.

The main difference between the Euclidean and the hyperbolic case is
the fact that
\begin{equation}
\label{A3} A_3^1 = A_3^2 = e^{\mp\frac{c_0^2}{2} \pi},
\end{equation}
i.e., $A_3^j$ grow exponentially with $c_0$ in the hyperbolic case.
\begin{figure}[!ht]
\centering
\resizebox{9cm}{5cm}{\includegraphics{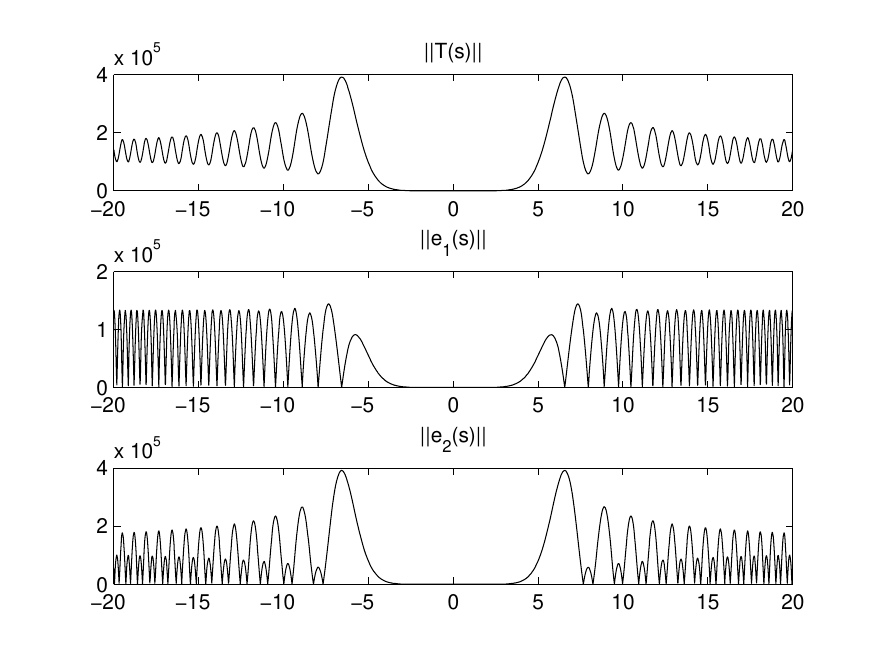}}
{\caption{\small $c_0 = 2.7$, hyperbolic case. Euclidean lengths of
$\T$, $\pp$ and $\qq$} \label{figure bound}}
\end{figure}
Thus, although $\T$, $\pp$ and $\qq$ are also bounded in the
hyperbolic case for a given $c_0$, this is by no means trivial,
unlike in the Euclidean case; moreover, there are no global bounds
valid for all $c_0$ and the bounds grow very fast even for
relatively small values of $c_0$. In Figure \ref{figure bound}, for
instance, we have plotted the Euclidean lengths of $\T$, $\pp$ and
$\qq$ in the hyperbolic case, for $c_0 = 2.7$; the lengths are seen
to be $\mathcal O(10^5)$. This fact makes the proof of Theorem
\ref{teorema} much more difficult in the Hyperbolic case, as well as
the numerical treatment of \eqref{xtposneg} and \eqref{ttposneg} for
large values of $c_0$.

Now, \eqref{xtposneg} is invariant under translations and also under
rotations in the Euclidean case and under Lorentz transformations
with unitarian determinant in the hyperbolic case; hence, given two
arbitrary $\A^j$ in $\mathbb S^2$ or $\HH$, we can transform
\eqref{xtposneg}, so that $\A^j$ have the form as in part (iii) of
Theorem \ref{teorema}. Then, it is possible to obtain the
corresponding $c_0$. Observe that for all $\A^j\in\HH$ there is a
corresponding $c_0$. However, when $\A^1=\A^2\in\mathbb S^2$, which
can be reduced to $\A^1=\A^2=(0, 0, 1)$, $c_0$ should be infinity,
so there is no real $c_0$ matching that case.

There is a natural connection between equations (\ref{xtposneg}) and
(\ref{ttposneg}), and the nonlinear Schr\"odinger equation. Equation
(\ref{xtposneg}) can be transformed into the Schr\"odinger equation
with a cubic nonlinearity via a transformation introduced by
Hasimoto \cite{hasimoto}. Especifically, define
\begin{equation}
\label{hasimoto-transform} \psi(s,t) = c(s,t) \exp {\left (i
\int_0^s \tau(s',t')\,ds' \right ) }.
\end{equation}
Then, $\psi$ satisfies the equation
\begin{equation}
\label{cubic-schrodinger} i \psi_t + \psi_{ss} \pm \frac{1}{2}\left
( |\psi|^2 +A(t)\right )\psi = 0,
\end{equation}
for some function $A(t)$, which can be removed by means of an
integrating factor.

Equation (\ref{ttposneg}) can be transformed into a nonlinear
Schr\"odinger equation by performing a stereographic projection onto
the complex plane. If we define
\begin{equation}
z(s,t) = \frac{T_1}{1+T_3} + i \frac{T_2}{1+T_3},
\end{equation}
then $z$ satisfies the equation \cite{sulem1}
\begin{equation}
\label{nonintegrable-schrodinger} z_t = iz_{ss} \mp
\frac{2i\bar{z}}{1\pm |z|^2}z_s^2.
\end{equation}
In this article, we study the self-similar solutions described above
from a numerical point of view, in order to understand the
singularity formation, as well as the mechanism for energy
concentration that is responsible for the formation of such
singularity. These solutions can be studied at the level of
$\mathbf{X}$ (eq. (\ref{xtposneg})), $\mathbf{T}$ (eq.
(\ref{ttposneg})), $\psi$ (eq. (\ref{cubic-schrodinger})), or at the
level of $z$ (eq. (\ref{nonintegrable-schrodinger})). However, we
will not work directly with (\ref{xtposneg}), given that
$\mathbf{X}$ can be recovered from $\mathbb R$  except for a
constant of integration that is fixed by considering \eqref{0xc0},
i.e.,
\begin{equation}
\X(s, 0) = 2c_0\sqrt t(0, 1, 0).
\end{equation}
Although there is a rich literature concerning the numerical study
of \eqref{ttpos} and, in general, the Landau-Lifshitz equation (see,
for instance, \cite{cimrak,carlos1,carlos2,wang}), the materials
analyzing the numerics of the self-similar solutions we are
considering are very scarce; the most relevant results were given by
Buttke some twenty years ago \cite{buttke1}. Nevertheless, Buttke
only considered the forward case, starting with a singular initial
datum at $t = 0$ and going forward in time; we are interested in the
singularity formation process, and therefore will focus mostly in
the backward case, i.e., starting from $t = 1$, we will approach the
singularity time.

Moreover, we are not aware of any other study of the  hyperbolic
case. For small $c_0$, the results of both cases are virtually
identical; because of that, we have used in our experiments $c_0 =
0.2$. In contrast, for bigger $c_0$, because of the rapid growth of
the Euclidean lengths of $\T$, $\pp$ and $\qq$ (see Figure
\ref{figure bound}), the hyperbolic case becomes much more difficult
to treat numerically, requiring a finer study, which we postpone for
the future.

The three equations (\ref{ttposneg}), (\ref{cubic-schrodinger}), and
(\ref{nonintegrable-schrodinger}) are posed in the whole real line. To
carry out the numerical simulations, we must restrict ourselves to an
interval $[-L,L]$ (for $L$ large enough). It is important then to give
appropriate boundary conditions in order to capture the singularity
formation mechanism. In \cite{buttke1}, only periodic boundary
conditions were considered.

The experiments performed here constitute numerical evidence of the
stability of \eqref{xtposneg}, \eqref{ttposneg} and
\eqref{nonintegrable-schrodinger}, and in particular of the robustness
of the energy concentration mechanism. Indeed, one can never introduce
numerically completely exact boundary conditions at $s = \pm L$; in the
best of the cases, there will be small errors due to machine precision
and in the worst one, we will be using very rough boundary conditions.
Therefore, we can consider that we are computing perturbed versions of
the exact solutions. In a recent series of articles
\cite{banica1,banica2} Banica and Vega prove the stability of the
solutions that we consider at the level of $\mathbf{T}$. At the level
of $\psi$, Banica and Vega \cite{banica2} also prove that a logarithmic
instability appears when $t$ comes close to zero. Notice that
$\mathbf{T}$ involves one integral with respect to the curvature and
two with respect to the torsion. In our case these are highly
oscillatory integrals so that, due to the cancelations, $\mathbf{T}$
and $\mathbf{X}$ become stable. These results give the theoretical
basis to the numerical results presented here, and suggest that it is
better to work numerically with the equations at the level of
$\mathbf{T}$ or $\mathbf{z}$. Therefore, we consider only equations
(\ref{ttposneg}), and (\ref{nonintegrable-schrodinger}).

We return now to the issue of boundary conditions. For equation
(\ref{ttposneg}), we consider two types of boundary conditions,
derived from the asymptotics of $\T$ given in Theorem \ref{teorema}.
The first boundary condition is simply $\T(\pm L, t) = \T(\pm, 1)$
and for the second one, we keep the first-order term in the
asymptotics of $\T(s, 1)$ and use the fact that $\T(s, t) =
\T(s/\sqrt t, 1)$.

For equation (\ref{nonintegrable-schrodinger}), we consider three
types of boundary conditions: the first one is the projected
boundary condition used for $\mathbf{T}$, The second boundary
condition is derived from the self-similarity condition, i.e., $z(s,
t) = z(s/\sqrt t, 1)$. Unlike the former one, this condition is not
an approximate one, so it introduces no noise, giving very clean
results. Finally, we derive a radiation boundary condition. This
condition seems to perform best, and in contrast to the previous two
conditions, can be used in the progressive case.

The remainder of the article is organized as follows: Equation
(\ref{ttposneg}) is considered in section \ref{section findif},
where we use a simple spatial discretization using finite
differences. The nonlinearity in the equation and the constraint
$\mathbf{T}$=1, make it difficult to construct efficient implicit
methods for these equations \cite{carlos1}. For the two boundary
conditions mentioned earlier,  we observe an stability constraint of
$|\Delta t| = \mathcal O(\Delta s^2)$. The first condition only
reproduces the singularity formation from a qualitative point of
view, while the second one gives good results for all $s$, even for
very small $t$, although a very small $\Delta s$ may be needed.

In section \ref{section spectral}, we consider equation
(\ref{nonintegrable-schrodinger}). This equation adapts well to a
pseudo-spectral method using Chebyshev polynomials; comprehensive
information about this family of polynomials can be found, for
instance, in \cite{fornberg} and \cite{gottlieb}. We use also a
Chebyshev point distribution, which is now much better suited, since
a denser concentration of points is required near the boundary for
small $t$ and, specially, for big $L$, which is when the
pseudo-spectral method really excels.

In order to measure the accuracy of our results, we compute the
error in the curvature obtained from $\T$ (which, when necessary,
can be recovered from $z$ as well), since the correct curvature has
the known value $c(s, t) = \frac{c_0}{\sqrt t}$.

Finally, the stability results of Banica and Vega are described in
section \ref{section stability}.
\section{Schr\"odinger map}
\label{section findif} We consider equation
\begin{align}
\label{ttposnegL}
\begin{cases}
\Tt(s, t) = \T(s, t)\propn\Tss(s, t),
\quad s \in[-L,+L] \\
\T(s, 1) = \T^0(s).
\end{cases}
\end{align}
In his Ph.~D. Thesis, Buttke considered the following
Crank-Nicholson type numerical scheme to integrate the Euclidean
version of \eqref{ttposnegL} \cite{buttke2}:
\begin{equation}
\label{buttke}
\frac{\T(s, t + \Delta t) - \T(s, t)}{\Delta t} = \frac{\T(s, t +
\Delta t) + \T(s, t)}{2} \wedge \left ( \frac{D_{+-} \T(s, t +
\Delta t) +D_{+-} \T(s,t)}{2} \right ),
\end{equation}
where $D_{+-}$ represents the approximation to the second derivative
using standard centered differences, i.e.,
\begin{equation}
D_{+-} T(s,t) = \frac{T(s+\Delta s,t)-2T(s,t)+T(s-\Delta
s,t)}{\Delta s^2}.
\end{equation}
Buttke only studied the forward case, i.e., starting from $t = 0$
and, concentrating all the information of $\T$ at $s = 0$, he tried
to recover the self-similar solutions we have described above, by
imposing periodic boundary conditions at $\pm L$.

The scheme used by Buttke's is implicit, and he used a fixed-point
iteration method to advance to the next time step. For this
iteration to converge, a time step $\Delta t = \mathcal O(\Delta
s^2)$ is needed. Therefore, although the method is a priori
unconditionally stable, in practice it is not.

We have observed that an explicit finite difference scheme is
equally efficient; more precisely, a scheme using a second-order
finite difference scheme in space with the classical fourth order
Runge-Kutta in time works well. In what follows, we have considered
the backward case, starting at $t = 1$ and trying to arrive at $t =
0$.

We divide $[-L,+L]$ in $N$ equally spaced parts, $ -L = s_0 < s_1 <
\cdots < s_{N-1} < s_N = +L$, with
\[ s_i = -L + i\Delta s, \qquad \Delta s = \frac{2L}{N}, \qquad i
= 0, \cdots, N.
\]
Equation \eqref{ttposneg} is discretized as
\begin{align*}
\Tt(s, t) & = \T(s, t)\propn D_{+-}T(s,t),
\end{align*}
and the time-stepping is carried out with the classical fourth order
Runge-Kutta. By using a fourth order Runge-Kutta, we guarantee that
$\T\cirpn\T = \pm 1$ is preserved with high accuracy. Nevertheless,
we normalize $\T$ at every time step by doing
\begin{align}
\label{renormalizarT} \T_i^{n + 1} & \equiv \dfrac{\T_i^{n + 1}}
{(\pm\T_i^{n + 1}\cirpn\T^{n + 1})^{\frac{1}{2}}}.
\end{align}
A higher order finite difference scheme could be effortlessly
implemented, although we have observed that this does not improve
significantly the quality of our results.

Experimentally, $|\Delta t| \lesssim 0.7\Delta s^2$ is found to be
needed for the method to be stable. Thus, as in Buttke's method, we
have a $\Delta t = \mathcal O(\Delta s^2)$ restriction, but it is
straightforward to advance to the next time step.

From theorem \ref{teorema}, we know that the self-similar solution
$\T^0$ has the following asymptotic expansion:
\begin{equation}
\label{Tasint} \T^0(s) = \mathbf A^j(c_0) - 2c_0\dfrac{\qq}{s} +
\mathcal O(1 / s^2).
\end{equation}
We approximate the boundary conditions for \eqref{ttposneg} by the
values of $T^0(s/\sqrt{t})$, and consider two types of boundary
conditions:
\begin{enumerate}
\item First order boundary conditions:
\begin{eqnarray}
T(+L,t) &=& A^{+}, \nonumber \\ T(-L,t) &=& A^{-}.
\end{eqnarray}
\item Second order boundary conditions:
\begin{eqnarray}
\label{second-order}
T(+L,t) &=& A^{+} - 2c_0 \sqrt{t} \frac{\qq(+L/\sqrt{t})}{L}, \nonumber \\
T(-L,t) &=& A^{-} + 2c_0 \sqrt{t} \frac{\qq(-L/\sqrt{t})}{L},
\end{eqnarray}
normalized to be of length $1$. In \eqref{second-order}, we use the
leading order in the asymptotics for $\qq$, obtained from theorem
\ref{teorema} on page \ref{teorema} as well.
\end{enumerate}

\subsection{First order boundary conditions}
We set
\begin{eqnarray}
T(+L,t) &=& A^{+}, \nonumber \\ T(-L,t) &=& A^{-},\ \ \ \forall t,
\end{eqnarray}
and obtain the initial datum $\T(s, 1)$ by integrating \eqref{Tpqt}
for $s\in[-L, L]$, using \eqref{Tpqt0} as initial condition. We have
also used a fourth order Runge-Kutta, normalizing every $\T^n \equiv
\T(s_n, 1)$, $\pp^n\equiv\pp(s_n, 1)$ and $\qq^n\equiv\qq(s_n, 1)$:
\begin{align}
\T^{n + 1} & \equiv \dfrac{\T^{n + 1}}{(\pm\T^{n + 1}\cirpn\T^{n +
1})^{\frac{1}{2}}}, & \pp^{n + 1} & \equiv \dfrac{\pp^{n +
1}}{(\pp^{n + 1}\cirpn\pp^{n + 1})^{\frac{1}{2}}}, & \qq^{n + 1}
\equiv \dfrac{\qq^{n + 1}}{(\pm\qq^{n + 1}\cirpn\qq^{n +
1})^{\frac{1}{2}}}.
\end{align}
We have executed our method $c_0 = 0.2$ and different value for $L$.
In our experiments, we have chosen $|\Delta t| = 0.5\Delta s^2$.
Smaller $|\Delta t|$ do not improve the results, since the method is
fourth-order accurate in time. With fixed boundary conditions, the
choice of $\Delta s$ is also not very important; in what follows, we
show the results for $L = 10$ and $L = 50$, having used $\Delta s =
0.01$, $\Delta t = -5\cdot10^{-5}$. With these parameters, we have
seen that smaller $\Delta s$ produce virtually identical results.
Nevertheless, for bigger $L$, it is convenient to use smaller
$\Delta s$.

To measure the quality of the results at a given $t$, we analyze the
curvature $c = \sqrt{\Ts\cirpn\Ts}$.

\begin{figure}[!ht]
\centering
\resizebox{9cm}{5cm}{\includegraphics{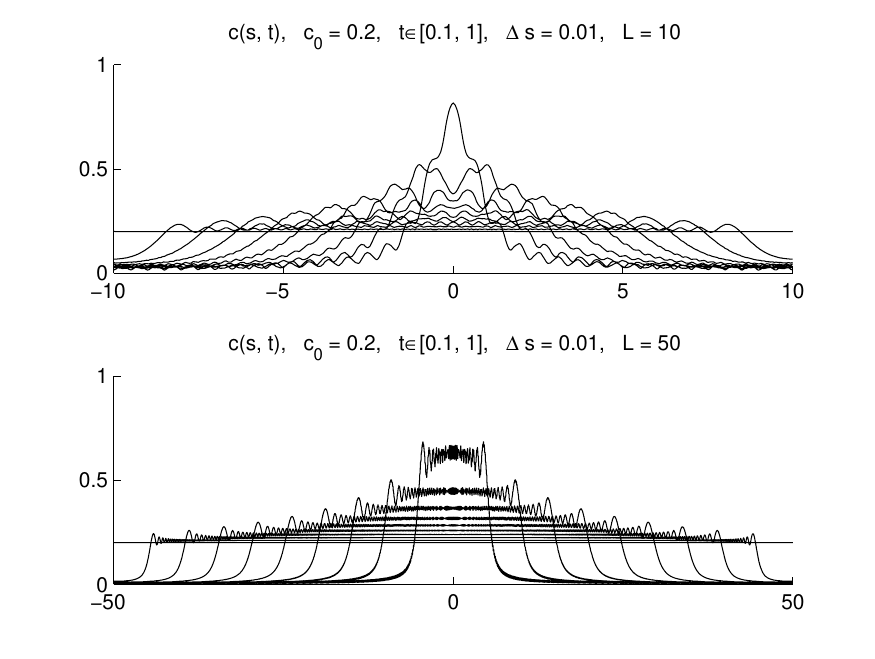}}
{\caption{\small Curvature at $t\in\{1, 0.9, \cdots, 0.1\}$, with
$\Delta s = 0.01$, $L = 10, 50$. Although the results for $L = 10$
are very poor, for $L = 50$ we notice a remarkable improvement.}
\label{figure curvature L}}
\end{figure}
In Figure \ref{figure curvature L} we show the curvature as a
function of $s$ and $t$, for $L = 10$ and $L = 50$. The accuracy is
lost as we approach $t = 0$, although as one would expect, the
bigger $L$ is, the better the results are. However, it is remarkable
the good accuracy with which one recovers the curvature at $s = 0$,
even for small times, as we can see in Figure \ref{figure curvature
s=0}.
\begin{figure}[!ht]
\centering
\resizebox{10cm}{5cm}{\includegraphics{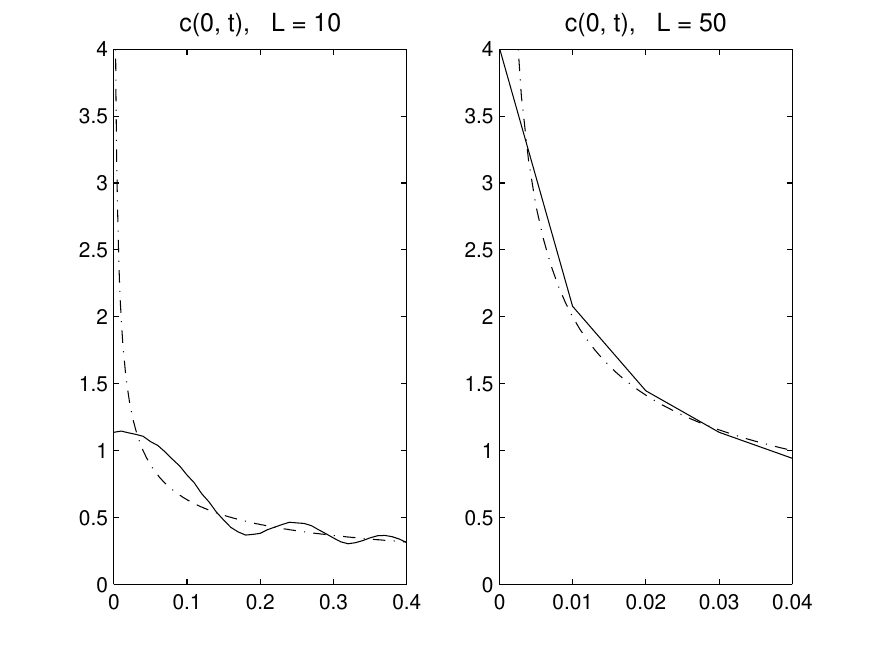}}
{\caption{\small Curvature at $s = 0$, with $\Delta s = 0.01$, $L =
10, 50$, together with its theoretical value. The accuracy improves
notably when increasing $L$.} \label{figure curvature s=0}}
\end{figure}
Again, this accuracy improves as we increase $L$. Finally, it is
also remarkable the fact that the energy between $[-L, L]$,
\begin{equation}
\int_{-L}^{+L}c^2(s, t)ds,
\end{equation}
calculated with the trapezoidal rule, is preserved with several
precision digits (see Figure \ref{figure error energy diffin
extremos fijos}). This is approximately $2Lc_0^2$ for all $t$.
\begin{figure}[!ht]
\centering
\resizebox{9cm}{5cm}{\includegraphics{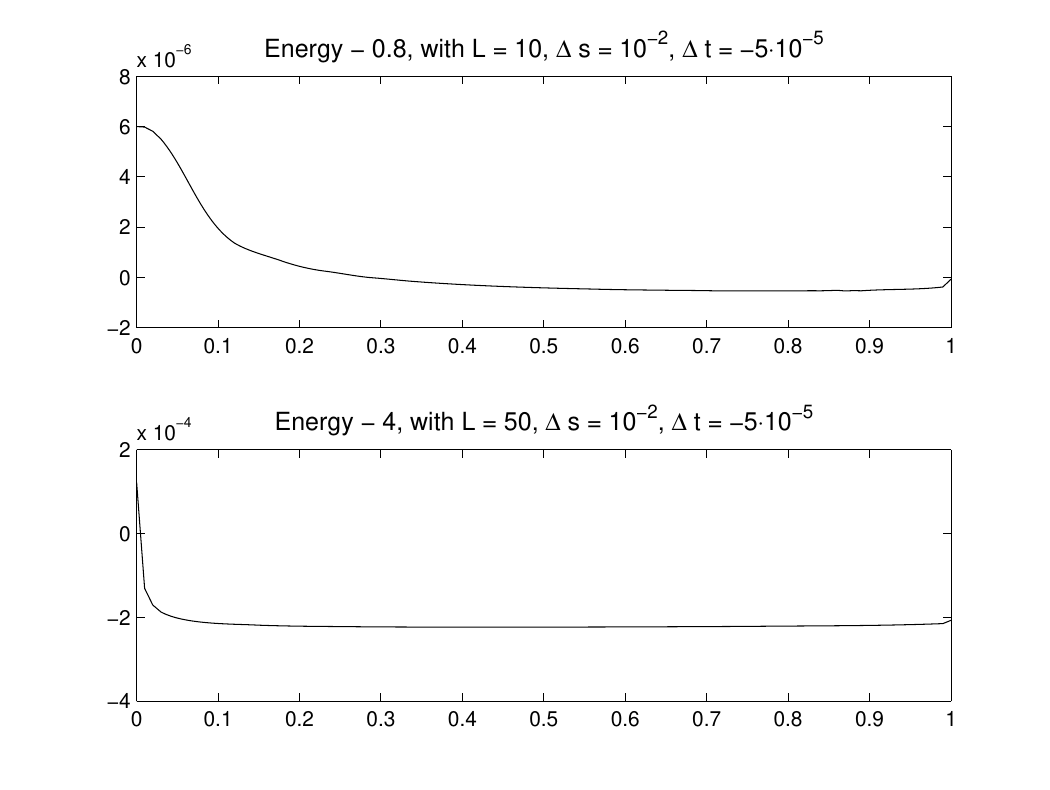}}
{\caption{\small Error in the conservation of the energy, for $L =
10$ and $L = 50$, using finite differences and fixed extremes.
Although the method is not symplectic, it preserves the energy very
well.} \label{figure error energy diffin extremos fijos}}
\end{figure}
Since the curvature at $s = 0$ is quite well recovered and the
energy well preserved, we are able to reproduce the behavior of $\T$
correctly, at least from a qualitative point of view. Indeed, in the
exact solution, the energy in $(-\infty, \infty)$ is preserved:
\begin{equation}
\int_{-\infty}^{+\infty}c^2(s, t)ds = \infty, \qquad \forall t;
\end{equation}

\noindent this infinite energy tends to concentrate at $s = 0$,
which causes the singularity to be produced. Now, we are working at
$[-L, L]$, but the finite energy is also preserved and it also tends
to concentrate on $s = 0$, so we can approximate the formation of
the singularity at $t = 0$. Moreover, by increasing $L$, we are able
to recover $c(0, t)$ for smaller times, improving the quality of the
results.

\subsection{Second order boundary conditions}
We set now
\begin{eqnarray}
\label{second-order-bc} \T(+L, t) &=& \A^+ - 2c_0\sqrt
t\dfrac{\mathbf e_2(L, t)}{L}
= \A^+ + 2c_0\sqrt t\dfrac{\Im [\B^+ e^{iL^2/4t}]}{L}, \\
\T(-L, t) &=& \A^- + 2c_0\sqrt t\dfrac{\mathbf e_2(-L, t)}{L} = \A^-
- 2c_0\sqrt t\dfrac{\Im [\B^- e^{iL^2/4t}]}{L},
\end{eqnarray}
normalized to be of length $1$.

Using these asymptotics as the boundary condition allows us to
obtain solutions that are not only qualitatively, but also
quantitatively correct. The reason for this is that the new boundary
conditions allow us to introduce energy into the system, as time
evolves, which concentrates in order to form a cusp at the origin
$s=0$.

\begin{figure}[!ht]
\centering
\resizebox{9cm}{5cm}{\includegraphics{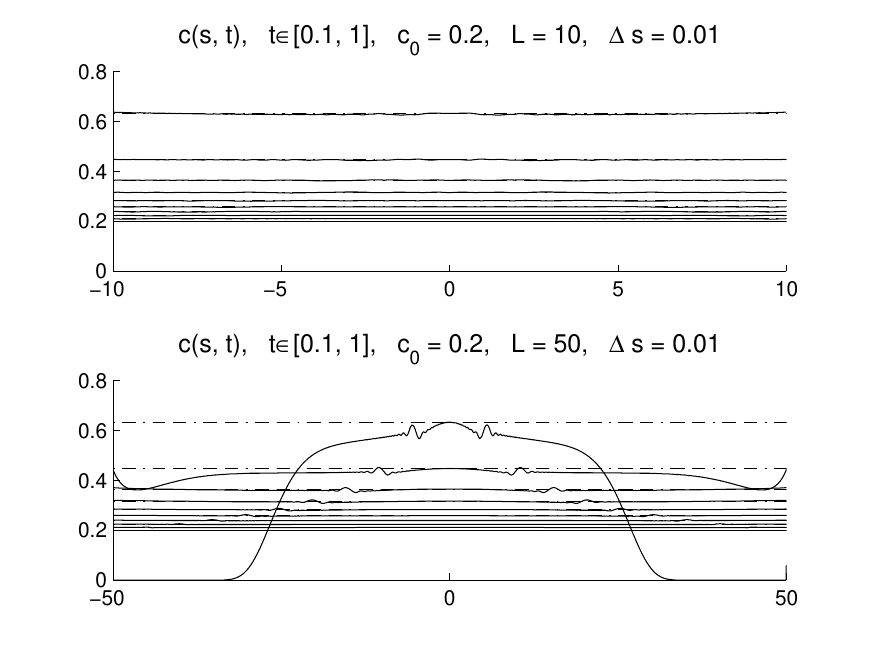}}
{\caption{\small Curvature in $t\in\{1, 0.9, \cdots, 0.1\}$, con
$\Delta s = 0.01$ y $L \in \{10, 50\}$. The results tend to worsen
when increasing $L$.} \label{figure curvdiffin L frontera}}
\end{figure}

Here, the choice of $\Delta s$ is much more delicate. For small $t$
and, specially, for big $L$, we have problems of aliasing, since the
boundary conditions \eqref{second-order-bc} are highly oscillatory.
In Figure \ref{figure curvdiffin L frontera}, we show the curvature
at several instants in time, obtained using $\Delta s = 0.01$,
$\Delta s = -5\cdot10^{-5}$ and $L = 10$ and $L = 50$. At $t = 0.1$,
the curvature has completely degraded for $L = 50$.

\subsection{Progressive case}

In the progressive case, considered by Buttke in his Ph.~D. thesis
\cite{buttke2}, we start at $t = 0$ and go forward in time, trying
to recover the self-similar solutions. Hence, we want to solve
numerically the following initial value problem
\begin{equation}
\label{Tt,t=0} \begin{cases}\Tt(s, t) = \T(s, t)\propn\Tss(s, t),
\quad s \in(-\infty,+\infty),
\\
\T(s, 0) = \A^1\chi_{[0,+\infty)}(s) + \A^2\chi_{(-\infty,0]}(s).
\end{cases}
\end{equation}

\noindent Again, since we cannot consider the whole $\mathbb R$, we
study
\begin{equation}
\begin{cases}
\Tt(s, t) = \T(s, t)\propn\Tss(s, t), \quad s \in[-L,+L],
\\
\T(s, 0) = \A^1\chi_{[0,+\infty)}(s) + \A^2\chi_{(-\infty,0]}(s).
\end{cases}
\end{equation}

\noindent We impose the boundary conditions $\T(+L, t) = \T(+L, 0) =
\A^1$ and $\T(-L, t) = \T(-L, 0) = \A^2$. Dividing $[-L, L]$ in $N$
parts of the same length, our initial datum will be
\begin{align}
\begin{cases}
\T_i^0 = \A^2, & i = 0, \cdots, \dfrac{N}{2} - 1, \\
\T_i^0 = (0, 0, 1), & i = \dfrac{N}{2}, \\
\T_i^0 = \A^1, & i = \dfrac{N}{2} + 1, \cdots, N,
\end{cases}
\end{align}

\noindent where $\T^0_i \equiv \T(s_i, t^0) = \T(s_i, 0)$, $s_i = -L
+ i\dfrac{2L}{N}$ and $\T(0, 0) = (0, 0, 1)$, in order to be
consistent with the backward case. Because of symmetries, $\T(0, t)
= (0, 0, 1)$, $\forall t > 0$. Because of part (iii) of Theorem
\ref{teorema}, $(\A^1)_3 = (\A^2)_3 = e^{\mp{\frac{c_0^2}{2}}\pi}$.
Therefore, once $c_0$ is fixed, we can chose some adequate $\A^1$
and $\A^2$:
\begin{align}
\label{2A+- en funcion de c_0} \A^1 & = \left(\sqrt{\pm(1 - e^{\mp
c_0^2\pi})}, 0, e^{\mp \frac{c_0^2}{2}\pi}\right), & \A^2 & =
\left(-\sqrt{\pm(1 - e^{\mp c_0^2\pi})}, 0, e^{\mp
\frac{c_0^2}{2}\pi}\right).
\end{align}

\noindent If we executed the backward case with that $c_0$, we would
obtain those $\A^1$ and $\A^2$, except for a rotation around the $z$
axis, which numerically has no relevance.

\begin{figure}[!ht]
\centering
\resizebox{9cm}{5cm}{\includegraphics{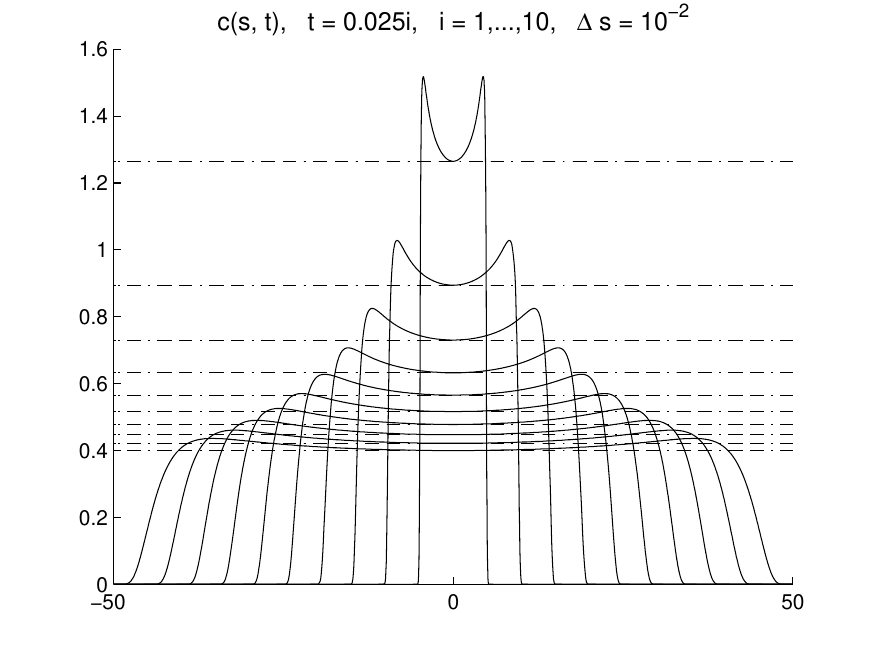}}
{\caption{\small Curvature for different $t$. Progressive case with
finite differences. We observe that the exactitude of the results
improves with time; the best result is obtained at $t = 0.25$.}
\label{figure diffin forward no conditions}}
\end{figure}

In Figure \ref{figure diffin forward no conditions}, we have taken
$L = 50$, $10000 + 1$ points, i.e., $\Delta s = 0.01$; $\Delta t =
5\cdot10^{-5}$ and $c_0 = 0.2$. We have drawn the curvature between
$t = 0.025$ and $t = 0.25$, with increments of $0.025$. It is clear
that the information goes outwards with constant velocity. The
correct theoretical curvature is shown with discontinuous strokes.
We can appreciate that the accuracy of the results improves with
time.

Thus, if we want to recover the curve for a given interval at a given
time, we have to consider a $\Delta s$ small enough and a $[-L, L]$ big
enough.

\begin{figure}[!ht]
\centering \resizebox{10cm}{5cm}{\includegraphics{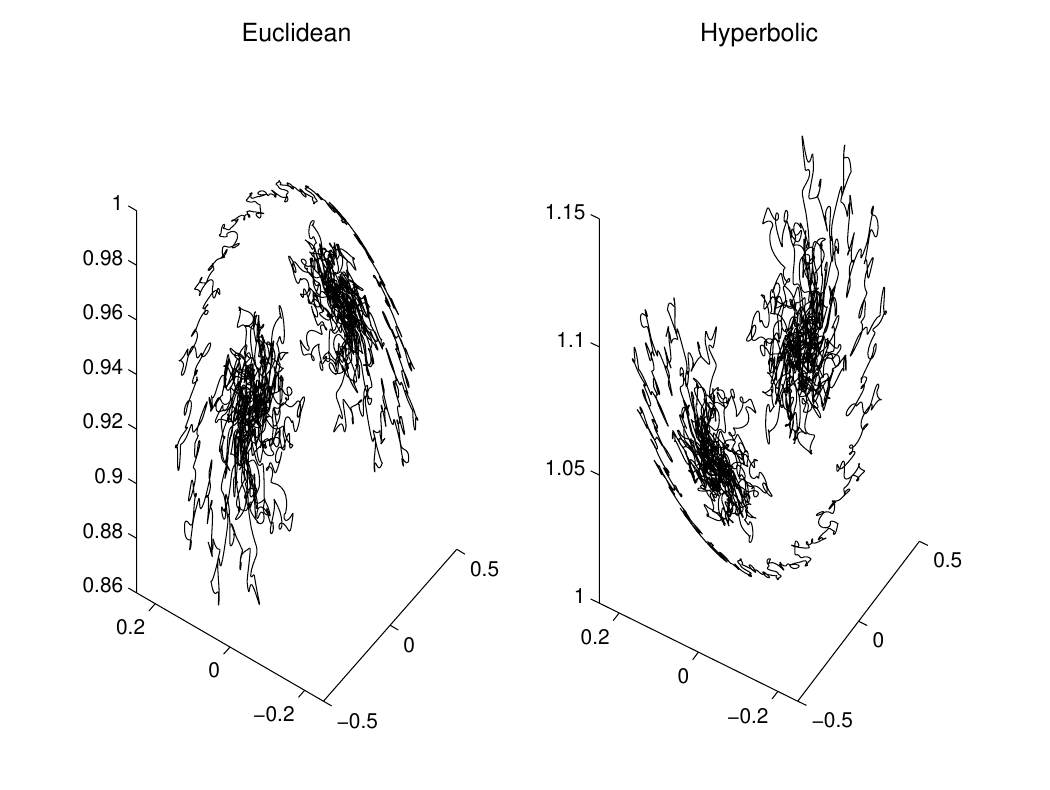}}%
\caption{Fractal formation in $\T$} \label{figure fractal T}
\end{figure}

If, however, we fix the boundary conditions, once the information
has reached the boundary, it is reflected back into the domain,
causing the appearance of fractal-like phenomena. In figure,
\ref{figure fractal T}, we have drawn the two cases of $\T$, for
$c_0 = 0.2$, $s\in[-50, 50]$, $\Delta s = 0.1$, $\Delta t =
5\cdot10^{-5}$, $t = 10$.
Although the approximated boundary
condition is convenient for the backward case, it does not seem to be appropriate
for the forward case, as we can see in Figure \ref{figure diffin
forward}.
\begin{figure}[!ht]
\centering
\resizebox{9cm}{5cm}{\includegraphics{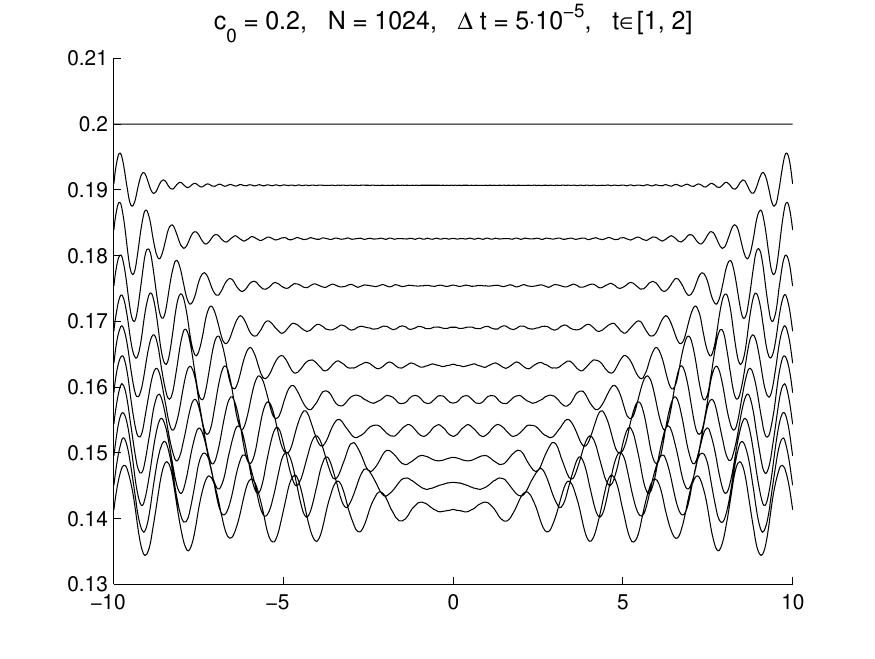}}
{\caption{\small Progressive case, starting from $t = 1$, till $t =
2$. As we see, the approximated boundary condition is not adequate
for the forward case.} \label{figure diffin forward}}
\end{figure}
We will get back to this issue in section \ref{section spectral},
where we consider a spectral discretization of
\eqref{nonintegrable-schrodinger}.

Finally, let us mention that the finite difference models have
another problem in the forward model, which Buttke also observed.
Indeed, the curvature error grows severely as we make $\Delta t\to
0$ (see Figure \ref{figure T noise curv}). In Section \ref{section
spectral}, this will also be solved.
\begin{figure}[!ht]
\centering \resizebox{9cm}{5cm}{\includegraphics{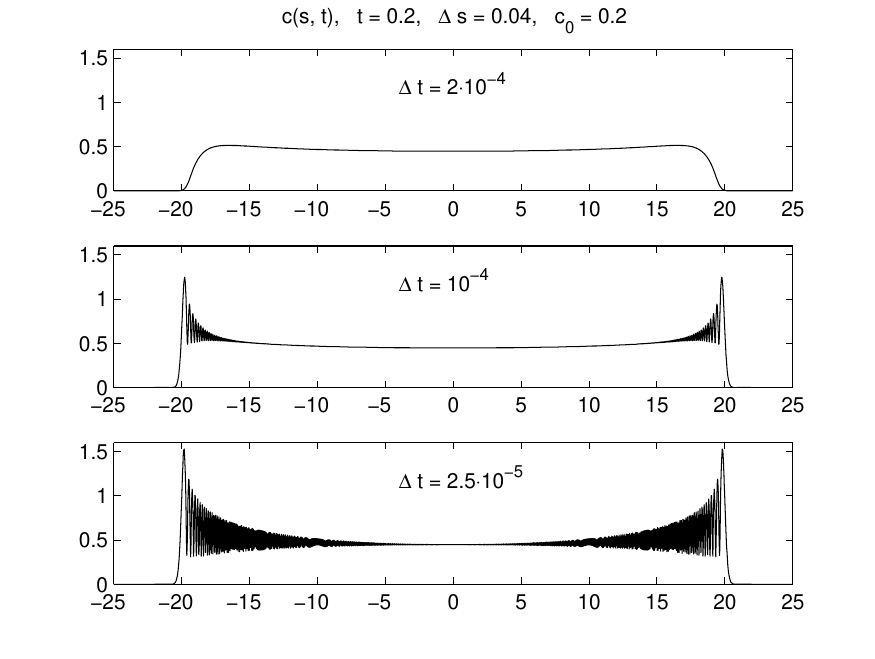}}%
\caption{Curvature noise when decreasing $\Delta t$} \label{figure T
noise curv}
\end{figure}

\subsection{Evolution of $\X$}

Once $\T(s, t)$ has being calculated, it is immediate to recover
$\X(s, t)$, since
\begin{align}
\begin{cases}
\Xs(s, t) = \T(s, t) \\
\X(0, t) = 2c_0\sqrt t(0, 1, 0);
\end{cases}
\end{align}
We integrate this equation using the scheme
\begin{align}
\label{2 integra X}\X^{n + 1} & = \X^n + \dfrac{\Delta s}{24}(9\T^n
+ 19\T^{n + 1} - 5\T^{n + 2} + \T^{n + 3}),
\\
\nonumber s_{n + 1} & = s_n + \Delta s,
\end{align}
which, in principle, is fourth order accurate in $\Delta s$. Note,
however, that $\mathbf{T}$ is only obtained to second order accuracy
in $\Delta s$, resulting overall in a second order method. However,
we use the fourth order scheme in order to reduce the accumulation
of errors.

\begin{figure}[!ht]
\centering
\resizebox{9cm}{5cm}{\includegraphics{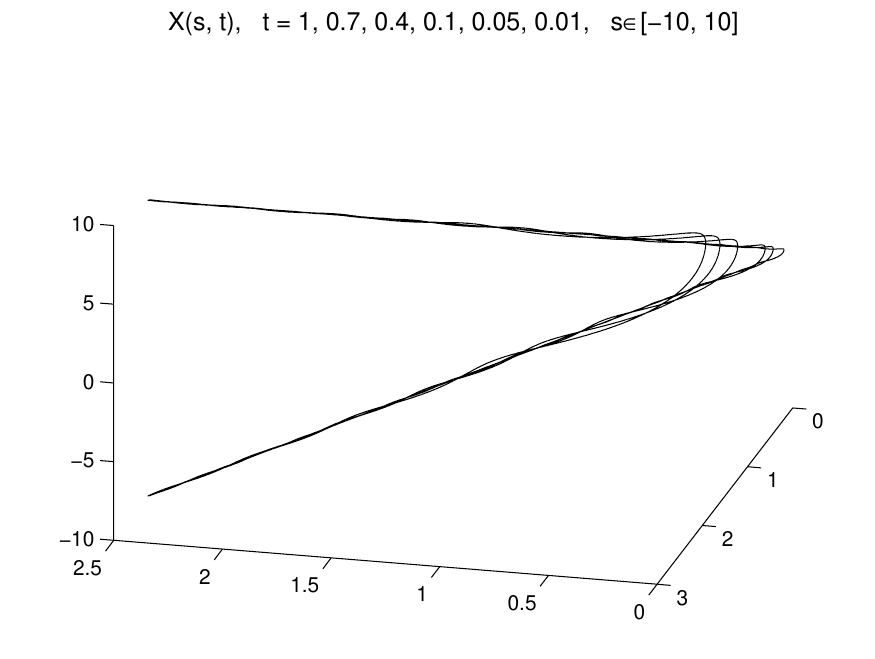}}%
\caption{\small Corner formation for $\X$. Fixed extremes.}
\label{figure X extremos fijos}
\end{figure}
We have recovered $\X$ for both boundary conditions. In Figure
\ref{figure X extremos fijos}, we have obtained $\T$ with fixed
extremes and the parameters of the experiment corresponding to the
second graphic in Figure \ref{figure curvature L}, i.e., Euclidean
case, $c_0 = 0.2$, $L = 50$, $\Delta s = 0.01$, $\Delta t =
-5\cdot10^{-5}$, drawing only $s\in[-10, 10]$.

In Figure \ref{figure X}, we have drawn the $\X$ recovered from
$\T$, with the parameters of the experiment corresponding to Figure
\ref{figure curvsTuperlast}. For this experiment, we had used the
second order boundary condition, so we can see the corner formation
with much greater clarity. Except for this, the two figures are
rather similar.
\begin{figure}[!ht]
\centering
\resizebox{9cm}{5cm}{\includegraphics{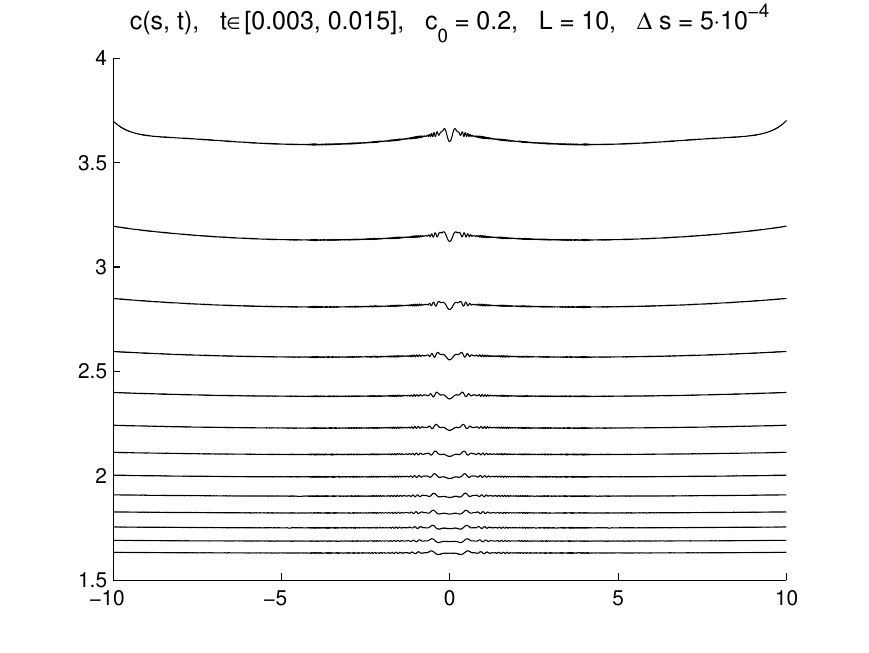}}
\caption{\small Curvature at $t \in\{0.003, 0.004, \cdots,
0.015\}$, with $\Delta s = 5\cdot10^{-4}$.} \label{figure
curvsTuperlast}
\end{figure}

\begin{figure}[!ht]
\centering \resizebox{9cm}{5cm}{\includegraphics{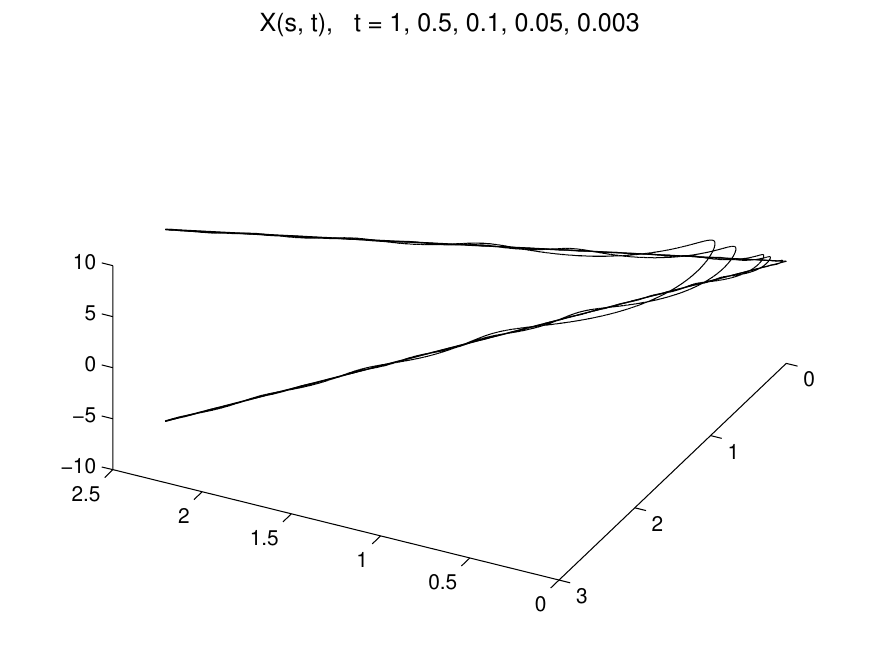}}%
\caption{Corner formation for $\X$. Approximated boundary.}
\label{figure X}
\end{figure}

\section{Schr\"odinger map: Stereographic projection}

\label{section spectral}

We consider now the stereographic projection of $\T=(T_1,T_2,T_3)$
over $\mathbb C$:
\begin{align}
\label{proy estereo} z = x + iy\equiv (x, y) \equiv
\left({\frac{T_1}{1 + T_3}}, {\frac{T_2}{1 + T_3}}\right).
\end{align}
We are projecting $\T$ from $(0, 0, -1)$ into $\mathbb R^2$,
identifying $\mathbb R^2$ with $\mathbb C$. In the Euclidean case,
where $\T\in\mathbb S^2$, there is a point on the sphere, $(0, 0,
-1)$, to which no point in $\mathbb C$ corresponds, because the
sphere is compact. Thus, we have a bijection between $\mathbb
S^2-\{(0, 0, -1)\}$ and $\mathbb R^2$. In the hyperbolic case, when
$\T\in\HH$, since $T_3
> 0$, we have a bijection between $\mathbb D$ and $\HH$, where
$$ \mathbb D = \{(x, y)\in\mathbb R^2\ /\ x^2 + y^2 < 1\} $$

\noindent is commonly referred to as the Poincar\'e disc. Unlike in
the Euclidean case, we do not need to eliminate any point from
$\HH$. The tangent vector $\T$ can be recovered from $z$ by the
inverse map:
\begin{equation}
\label{proy invers} \T = (T_1, T_2, T_3) \equiv \left({\frac{2x}{1
\pm x^2 \pm y^2}}, {\frac{2y}{1 \pm x^2 \pm y^2}}, {\frac{1 \mp x^2
\mp y^2}{1 \pm x^2 \pm y^2}}\right).
\end{equation}
Now, differentiating $z$ in \eqref{proy estereo}, together with
\eqref{ttposneg} and \eqref{proy invers}, we get the following
nonintegrable, nonlinear Schr\"odinger equation for $z$:
\begin{align}
\label{z_t} z_t = iz_{ss} \mp \frac{2i\bar z}{1 \pm |z|^2}z_s^2.
\end{align}
The advantage of this equation, as opposed to \eqref{xtposneg} and
\eqref{ttposneg}, is that the higher order term $z_{ss}$ treated
implicitly, eliminating,  or at least reducing significantly  the
restrictions for $\Delta t$.

We can also  express the other elements of the generalized
Frenet-Serret formulae \eqref{frenet-serretgen}, i.e., $c$, $\tau$,
$\pp$ and $\qq$ in terms of $z$.  We get the generalized curvature
from $c = (\Ts\cirpn\Ts)^{1/2}$:
\begin{align}
\label{cT} c = (T_{1s}^2 + T_{2s}^2 \pm T_{3s}^2)^{\frac{1}{2}} =
\frac{2(x_s^2 + y_s^2)^{\frac{1}{2}}}{1 \pm x^2 \pm y^2} =
\frac{2|z_s|}{1 \pm |z|^2}.
\end{align}
The expression for the generalized torsion is somewhat more
involved:
\begin{align}
\label{tauT} \nonumber\tau & =
\frac{1}{c^2}\T\circ_\pm(\Ts\propn\Tss)
\\
& = \frac{2\Im(z\bar z_s)}{|z|^2 \pm 1} + \frac{\Im(\bar
z_sz_{ss})}{ |z_s|^2}.
\end{align}
The expression for $\pp$ follows from $\pp = \frac{1}{c}\Ts$:
\begin{equation}
\label{ppz} \pp = \frac{1}{(x_s^2 +
y_s^2)^{\frac{1}{2}}}\left(\frac{x_s(1 \mp x^2 \pm y^2) \mp
2xyy_s}{1 \pm x^2 \pm y^2}, \frac{y_s(1 \pm x^2 \mp y^2) \mp
2xyx_s}{1 \pm x^2 \pm y^2}, \frac{\mp 2xx_s \mp 2yy_s}{1 \pm x^2 \pm
y^2}\right).
\end{equation}
Finally, we obtain $\qq$ from $\qq = \T\propn\pp$:
\begin{equation}
\label{qqz} \qq = \frac{1}{(x_s^2 + y_s^2)^{1/2}}\bigg(\frac{\mp
2xx_sy - y_s(1 \mp x^2 \pm y^2)}{1 \pm x^2 \pm y^2}, \frac{\pm
2xyy_s + x_s(1 \pm x^2 \mp y^2)}{1 \pm x^2 \pm y^2}, \frac{\pm 2xy_s
\mp 2x_sy}{1 \pm x^2 \pm y^2}\bigg).
\end{equation}

\subsection{Self-similar solutions}

\label{subsection self-similar}

The self-similar solutions of $\T$ correspond to self-similar
solutions of
\begin{align}
\label{z_t2} z_t = iz_{ss} \mp \frac{2i\bar z}{1 \pm |z|^2}z_s^2.
\end{align}
To obtain them, note that if $z$ is a solution of \eqref{z_t2}, so
is $z_\lambda(s, t) = z(\lambda s, \lambda^2 t)$, for all $\lambda$.
Taking $\lambda = t^{-1/2}$,
\begin{align*}
z(s, t) = z(s / \sqrt t, 1) = f(s /\sqrt t),
\end{align*}

\noindent where $f(s)\equiv z(s, 1)$. Introducing $f(s / \sqrt t)$
in \eqref{z_t2},
\begin{align*}
-\frac{1}{2}st^{-3/2}f'(s / \sqrt t) = \frac{i}{t}f''(s / \sqrt t)
\mp \frac{2i\bar f(s / \sqrt t)}{1 \pm |f(s / \sqrt
t)|^2}\frac{1}{t} (f'(s / \sqrt t))^2.
\end{align*}
Setting $t = 1$,
\begin{align*}
-\frac{1}{2}sf'(s) = if''(s) \mp \frac{2i\bar f(s)}{1 \pm |f(s)|^2}
(f'(s))^2;
\end{align*}

\noindent hence
\begin{align}
\label{z(s,1)} f''(s) = \frac{is}{2}f'(s) \pm \frac{2\bar f(s)}{1
\pm |f(s)|^2} (f'(s))^2.
\end{align}
Using the same initial conditions as before,  $\T(0, 1) = (0, 0, 1)$
and $\T_s(0, 1) = c_0(1, 0, 0)$, we get
\begin{align}
\begin{cases}
f(0) = z(0, 1) = 0, \\
f'(0) = z_s(0, 1) = \dfrac{c_0}{2}.
\end{cases}
\end{align}

\noindent It is immediate to generalize these expressions for times
other than $t = 1$. Defining $g(s) = f(s / \sqrt t) = z(s, t)$, the EDO
for $g(s)$ is
\begin{align}
\label{z(s,t)}
\begin{cases}
g''(s) = i\dfrac{s}{2t}g'(s) \pm \frac{2\bar g(s)}{1 \pm |g(s)|^2}
(g'(s))^2
\\
g(0) = f(0) = 0 \\
g'(0) = \dfrac{1}{\sqrt t}f'(0) = \frac{c_0}{2\sqrt t}.
\end{cases}
\end{align}

\noindent Integrating \eqref{z(s,1)}, or eventually \eqref{z(s,t)},
we get the initial datum for \eqref{z_t2}. We use a fourth-order
Runge-Kutta method to integrate \eqref{z(s,t)}. We obtain $f(s)$ for
$s \ge 0$. The values $f(s)$ for $s<0$ are obtained by symmetry,
since $f(s)$ is antisymmetric.

\subsection{A spectral collocation method for $z$}
We consider a semi-implicit method for
\begin{align}
\label{z_t3} z_t = iz_{ss} \mp \frac{2i\bar z}{1 \pm |z|^2}z_s^2,
\end{align}
where we treat the linear term on the right-hand side, $iz_ss$,
implicitly, and the nonlinear term explicitly. We use a Chebyshev
spectral collocation method \cite{Canuto:88}, with nodes $s_i$:
$$
s_i = L\cos\left(\frac{i\pi}{N}\right), \qquad i = 0, \cdots, N,
$$
and approximate $z$ by a polynomial of the form
$$
z(s, t) \approx \sum_{k = 0}^N a_k(t)T_k(s/L),
$$
where $T_x(s) = \cos( k \arccos s)$ is the Chebyshev polynomial of
degree $k$. The coefficients $\{a_k\}_{k=0}^N$ are obtained using
the fast Fourier transform (FFT) \cite{FFTW05}.

For the time evolution, we have chosen a second order, semi-implicit
Backward Differentiation Formula (BDF) \cite{lambert:73,ascher}:
\begin{equation}
\label{SBDF}\frac{1}{2\Delta t}\left[3U^{n + 1} - 4U^n + U^{n -
1}\right] = iU_{ss}^{n + 1} + 2\mathcal N(U^n, t^n) - \mathcal
N(U^{n - 1}, t^{n - 1}),
\end{equation}
where we denote by $\mathcal N(U,t)$ the nonlinear term on the
right-hand side of \eqref{z_t3}. The BDF is particularly suited for
this problem, as it imposes a very strong decay in the high
frequency modes.

Boundary conditions of Dirichlet or Neumann type can be easily
implemented:\begin{enumerate}
\item Dirichlet boundary condition:
\begin{equation}
\begin{cases}
u(-L, t) = u_1(t), \\ u(+L, t) = u_2(t)
\end{cases}
\longleftrightarrow
\begin{cases}
\sum_{k = 0}^N (-1)^k a_k^{n + 1} = u_1(t^{n + 1}) \\
\sum_{k = 0}^N a_k^{n + 1} = u_2(t^{n + 1}).
\end{cases}
\end{equation}
\item Neumann boundary condition:
\begin{equation}
\begin{cases}
u_s(-L, t) = u_1(t) \\ u_s(+L, t) = u_2(t)
\end{cases}
\longleftrightarrow
\begin{cases}
\dfrac{1}{L}\sum_{k = 0}^N (-1)^{k + 1}k^2a_k^{n + 1} = u_1(t^{n + 1}) \\
\dfrac{1}{L}\sum_{k = 0}^Nk^2a_k^{n + 1} = u_2(t^{n + 1}).
\end{cases}
\end{equation}
\end{enumerate}
We use a spectral filter, and set equal to zero all those coefficients
$a_k^{n + 1}$ whose modulus is smaller than a given $\varepsilon$
\cite{hou}. Specifically, we have made $a_k \equiv 0$ whenever $|a_k^{n
+ 1}| < 10^{-14}$.

For scheme \eqref{SBDF}  we need two initial conditions, one at time
$t = t^0$, which is known, and the other one at $t = t^1$. We obtain
$U(s, t^1)\approx u(s, t^1)$ using semi-implicit Backward Euler, and
Richardson extrapolation.

\subsection{Projected second order boundary condition}

In order to compare the finite-difference and the pseudo-spectral
method, we have projected the same second order boundary condition
we used for $\T$: We approximate $\T(L, t^{n + 1})$ as in
\eqref{second-order-bc}:
\begin{align*}
\tilde\T^{n + 1}(L) & = \tilde\A + 2c_0\sqrt{t^{n + 1}} \dfrac{\Im
[\tilde\B e^{iL^2/4t^{n + 1}}]}{L}, & \T^{n + 1}(L) & =
\dfrac{\tilde\T^{n + 1}(L)}{(\pm\tilde\T^{n + 1}(L)
\cirpn\tilde\T^{n + 1}(L))^\frac{1}{2}};
\end{align*}
and project it over $\mathbb C$, so
\begin{align*}
z^{n + 1}(L) = \dfrac{T^{n + 1}_1(L)}{1 + T^{n + 1}_3(L)} +
i\dfrac{T^{n + 1}_2(L)}{1 + T^{n + 1}_3(L)}.
\end{align*}

\noindent $\tilde\A$ and $\tilde\B$ are computed also using the
asymptotic expansion for $\T$:
\begin{align*}
\tilde\A & \equiv \T(L, t_0) + 2c_0 \sqrt t_0\dfrac{\mathbf e_2(L,
t_0)}{L},
\\
\tilde\B & \equiv (\pp(L, t) - i\qq(L, t))e^{-iL^2/4t},
\end{align*}

\noindent where $\T(L, t_0)$, $\pp(L, t_0)$ and $\qq(L, t_0)$ are
obtained from \eqref{proy invers}, \eqref{ppz} and \eqref{qqz},
respectively.

We choose the same parameters as in the finite difference case,
i.e.,  $c_0 = 0.2$. Only the results for the Euclidean case will be
given, since the results for the hyperbolic case are virtually
identical for this value of $c_0$. To measure the accuracy of our
results, we compute the curvature as a function of space and time,
and compare with the exact value $\frac{c_0}{\sqrt t}$. Thus, the
smaller the time that we reach with correct curvatures, the better
we will consider the results to be.
\begin{figure}[!ht]
\centering
\resizebox{9cm}{5cm}{\includegraphics{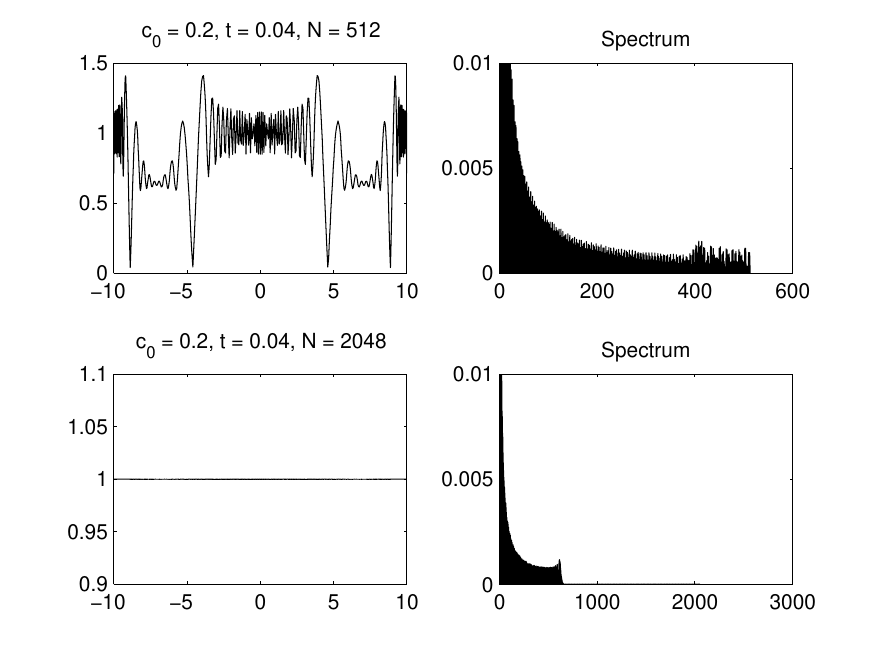}}
{\caption{\small Curvature and spectrum of the exact $z(s, 0.01)$,
when $N$ is too small ($N = 512$) and when $N$ is big enough ($N =
2048$).} \label{figure spectrum z}}
\end{figure}
To illustrate the need for high resolution, we show in figure
\ref{figure spectrum z} the curvature $c$ obtained for $t=0.04$. The
theoretical value is $c_0/\sqrt{t}=1$. As can be seen, when $N = 512$,
the accuracy is completely lost. In contrast, when $N = 2048$, the
error  is of the order of $10^{-8}$.

In all our numerical simulations, even for big $N$, we have found no
restriction for $|\Delta t|$, so there is evidence that the method
is unconditionally stable. Nevertheless, since the method is of
order two, diminishing $|\Delta t|$ can improve greatly the results,
provided $N$ is large enough. In figure \ref{figura 2 frontera
articulo dt = 1_5}, we plot the curvature for $N = 1024$ and $N =
2048$, at $t = 0.03$, $t = 0.04$ y $t = 0.05$. For $N = 1024$, the
noise indicates that the resolution is not high enough. Doubling the
number of collocation points ameliorates the situation.
\begin{figure}[!ht]
\centering
\resizebox{9cm}{5cm}{\includegraphics{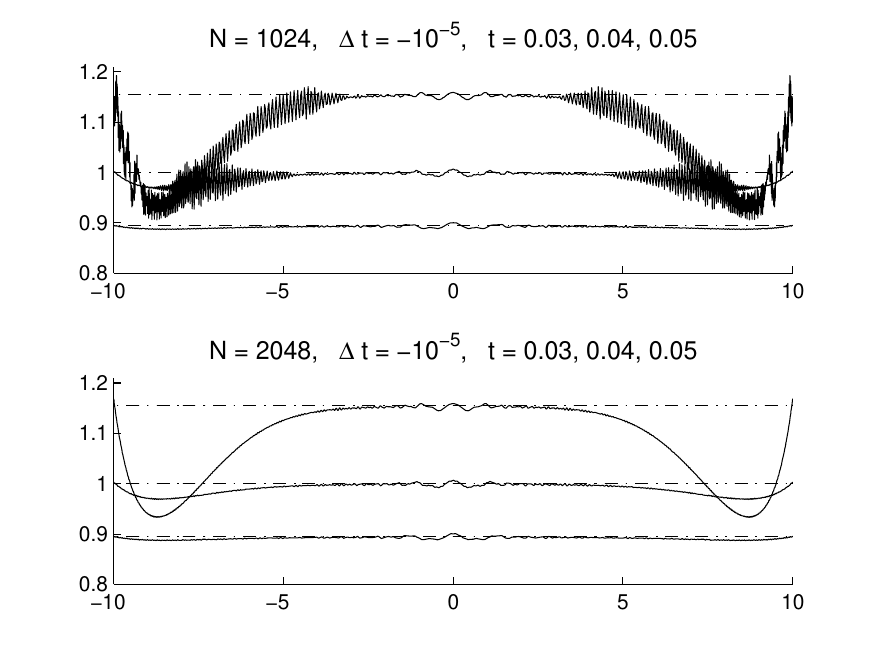}}
{\caption{\small Curvature, for $|\Delta t| = 10^{-5}$. When $N =
1024$, the spectrum has overflowed and there appears noise, because
of the lack of frequencies. With $N = 2048$, this noise have
disappeared.} \label{figura 2 frontera articulo dt = 1_5}}
\end{figure}
Reducing the time step is not enough to remove the oscillations:
using  $\Delta t = -10^{-6}$, the results at the same time instants
are practically exact for $N = 2048$, but, for $N = 1024$ the noise
has not disappeared, as we see in Figure \ref{figura 2 frontera
articulo dt = 1_6}.
\begin{figure}[!ht]
\centering
\resizebox{9cm}{5cm}{\includegraphics{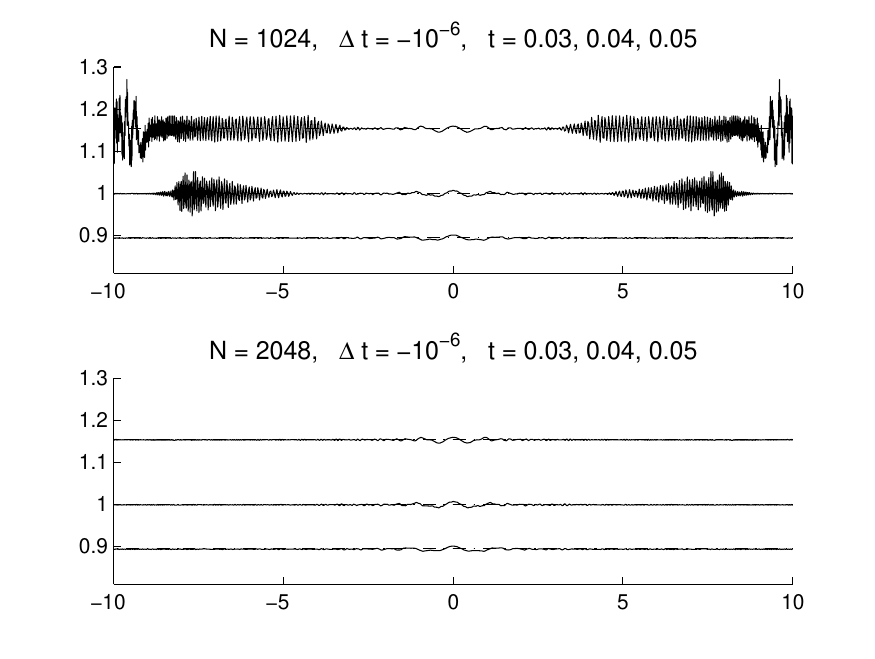}}
{\caption{\small Curvature, for $|\Delta t| = 10^{-6}$. Diminishing
$\Delta t$ is only useful when the number of frequencies is
satisfactory.} \label{figura 2 frontera articulo dt = 1_6}}
\end{figure}

In order to study the cusp-formation at $t=0$, we hade developed an
adaptive methodology, both in time and space. A simple strategy is to
simply duplicate the number of Chebyshev nodes whenever the derivative
of $z$ develops, high frequency components. Specifically, if we write
\[
z_s(s) = \sum_{k = 0}^Nb_kT_k(s / L),
\]
we duplicate the number of nodes when
\begin{align}
\max_{k\in[\frac{3N}{4}, N]}|b_k| > 2\cdot10^{-4}.
\end{align}
Whenever we duplicate the frequencies, we divide $\Delta t$ by $4$.

Using this strategy, we have solved the equation starting with
$N=1024$, and $\Delta t = -2 \times 10^{-6}$, up to $N=16384$.

Since our method is only second order in time, we would need a priori
small $|\Delta t|$ in order to decrease the numerical errors.
\begin{figure}[!ht]
\centering
\resizebox{9cm}{5cm}{\includegraphics{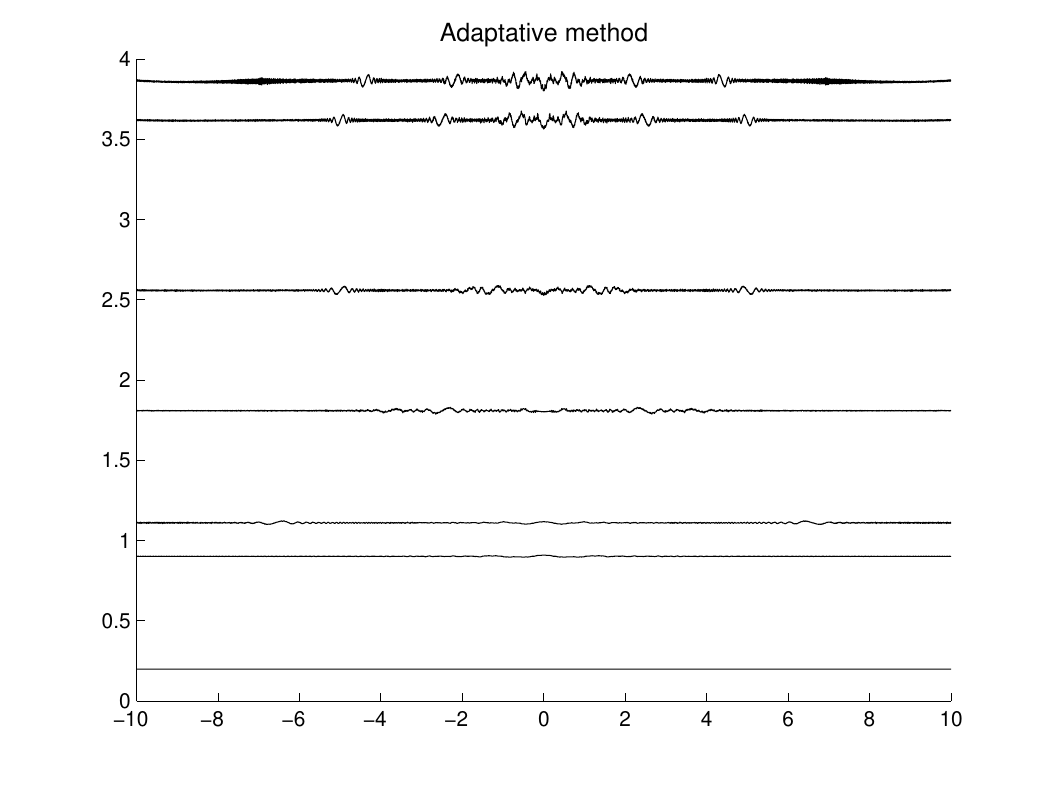}}
{\caption{\small Adaptive method, with approximated boundary
condition, from $N = 1024$ to $N = 16384$.} \label{figure frontera
articulo fulladapt}}
\end{figure}
In Figure \ref{figure frontera articulo fulladapt}, we plot the
curvature in those instants where we have duplicated number of
Chebyshev nodes, that is,  at $t = 4.91\cdot10^{-2}$, $t =
3.24\cdot10^{-2}$, $t = 1.22\cdot10^{-2}$ and $t =
6.09\cdot10^{-3}$, and gone down up to $t = 2.67\cdot10^{-3}$. Note
that this means that the curvature has been increased by a factor of
$20$, and the energy by a factor of $400$.

\subsubsection{Self-similarity boundary conditions}
\label{sec:self-similar-bc} Since we are trying to approximate the
self-similar solutions, i.e., $z(s, t) = z(s / \sqrt t, 1)$, it
seems natural to introduce this condition on the boundary. Upon
differentiation, we get 21\begin{align} z_t(s, t) =
-\dfrac{s}{2t}z_s(s, t),
\end{align}
which can be translated into the following boundary conditions:
\begin{align*}
\frac{z^{n + 1}(L) - z^{n - 1}(L)}{2\Delta t} & = -\frac{
L}{2t^n}z_s^n(L), \\
\frac{z^{n + 1}(-L) - z^{n - 1}(-L)}{2\Delta t} & = \frac{
L}{2t^n}z_s^n(-L).
\end{align*}
This choice of the boundary condition produces good results. However,
for accuracy reasons, for small $t$, we need to choose $\Delta t$
small enough to avoid numerical artifacts on the boundary. This is
illustrated in figure \ref{figure frontera autosemejante delta t
grande}, where we show what happens on the boundary when $\Delta t$
is chosen to be too large.
\begin{figure}[!ht] \centering
\resizebox{9cm}{5cm}{\includegraphics{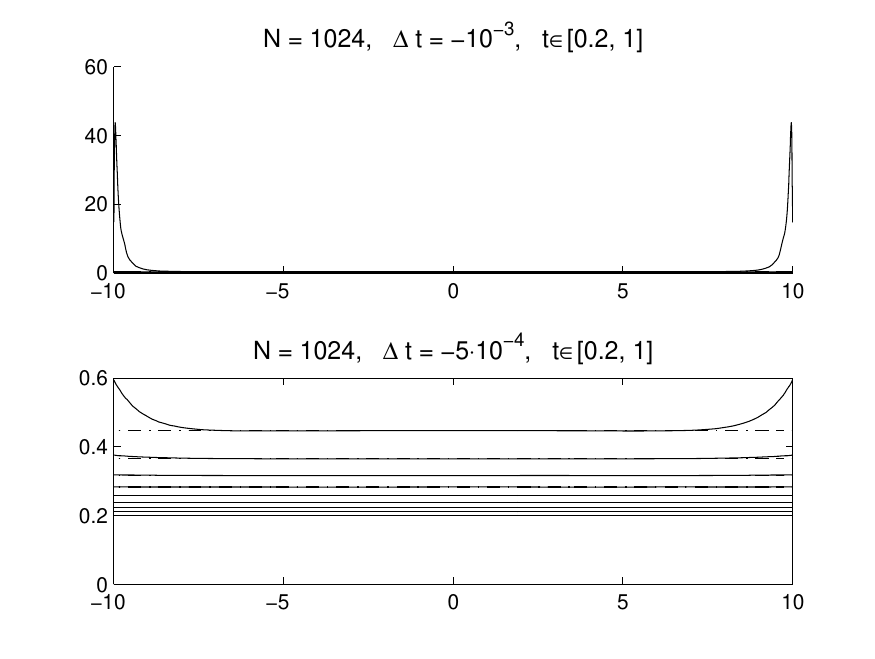}}
{\caption{\small Curvature for too big $\Delta t$. When $\Delta t =
-10^{-3}$, curvature explodes at $s = \pm L$ around $t = 0.2$. With
$\Delta t = -5\cdot10^{-4}$, there is a good improvement.}
\label{figure frontera autosemejante delta t grande}}
\end{figure}

\begin{figure}[!ht]
\centering
\resizebox{9cm}{5cm}{\includegraphics{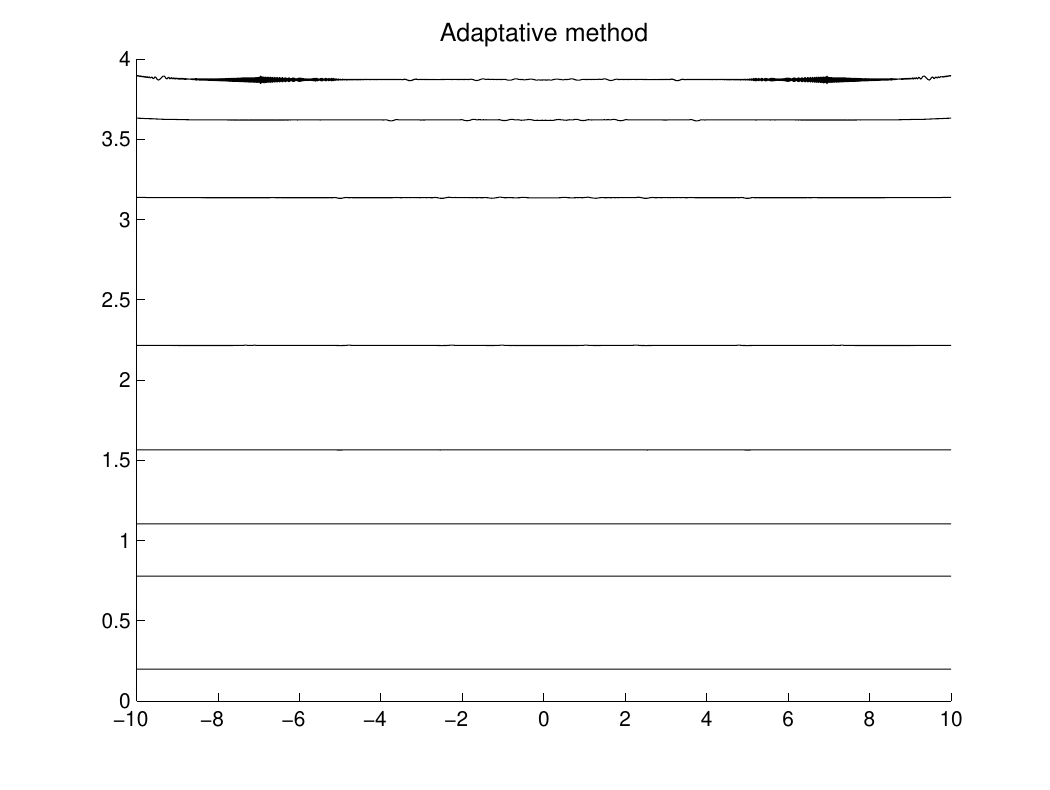}}
{\caption{\small Adaptive method, with self-similar boundary
conditions, from $N = 1024$ till $N = 16384$.} \label{figure
frontera autosem fulladapt}}
\end{figure}

In Figure \ref{figure frontera autosem fulladapt}, we have
implemented an adaptive method, starting from $N = 1024$, $\Delta t
= -2\cdot10^{-6}$. Representing $z_s(s) = \sum_{k = 0}^Nb_kT_k(s /
L)$, we have duplicated the frequencies when
\begin{align}
\label{2duplicaautosem} \max_{k\in[\frac{3N}{4}, N]}|b_k| >
5\cdot10^{-5};
\end{align}

\noindent i.e., at $t = 6.61\cdot10^{-2}$, $t = 3.23\cdot10^{-2}$, $t =
1.63\cdot10^{-2}$ and $t = 8.15\cdot10^{-3}$. With $N = 16384$, the
results are still valid for much smaller times; in the figure we have
also shown the curvature when $t = 4.07\cdot10^{-3}$, $t =
3.05\cdot10^{-3}$ and $t = 2.67\cdot10^{-3}$, multiplying almost by
$20$ the initial curvature. Observe that, in comparison with Figure
\ref{figure frontera articulo fulladapt}, the results of Figure
\ref{figure frontera autosem fulladapt} are much cleaner.

\subsection{Radiation boundary condition}
We derive another boundary condition, that aims at capturing the
correct flow of energy through the boundary at $s=\pm L$. For this,
we consider the Hasimoto transform, $\psi =
c\exp\left(i\int_0^s\tau(s', t)ds'\right)$. From \eqref{cT} and
\eqref{tauT}, it follows that
\begin{align*}
\label{psi1} \psi = \frac{2|z_s|}{1 \pm
|z|^2}\exp\left[i\int_0^s\left(\frac{2(yx_s - xy_s)}{\pm1 + x^2 +
y^2} + \frac{x_sy_{ss} - y_sx_{ss}}{x_s^2 + y_s^2}\right)ds'\right].
\end{align*}

\noindent Now,
\begin{align*}
\frac{x_sy_{ss} - y_sx_{ss}}{x_s^2 + y_s^2} = \frac{\displaystyle
\frac{y_{ss}}{x_s} - \frac{y_sx_{ss}}{x_s^2}}{1 +
\left(\displaystyle \frac{y_s}{x_s} \right)^2} =
\frac{\partial}{\partial s}\arctan\left(\frac{y_s}{x_s}\right),
\end{align*}

\noindent so
\begin{align}
\nonumber \psi & = \frac{2|z_s|}{1 \pm
|z|^2}\exp\left\{i\int_0^s\left[\frac{2(yx_s - xy_s)}{\pm1 + x^2 +
y^2} + \frac{\partial}{\partial s}
\arctan\left(\frac{y_s}{x_s}\right) \right]ds'\right\}
\\
\nonumber & = \frac{2|z_s|}{1 \pm
|z|^2}\exp\left[i\arctan\left(\frac{y_s}{x_s}\right) -
i\arctan\left(\frac{y_s(0)}{x_s(0)}\right)\right]
\exp\left[i\int_0^s\frac{2(yx_s - xy_s)}{\pm1 + x^2 + y^2}ds'\right]
\\
& = \frac{2z_s}{1 \pm |z|^2}\exp\left [i\int_0^s \frac{2(yx_s -
xy_s)}{\pm 1 + x^2 + y^2}ds'\right]\exp\left[-
i\arctan\left(\frac{y_s(0)}{ x_s(0)}\right)\right].
\end{align}
At $s=\pm L$, we know that $\psi(\pm L, t) = \dfrac{c_0}{\sqrt t}
e^{i\frac{L^2}{4t}}$. From this, it follows that
\begin{align*}
\frac{c_0}{\sqrt t}e^{i\frac{L^2}{4t}} = \left\{\frac{2z_s}{1 \pm
|z|^2}\exp\left[i\int_0^s\frac{2(yx_s - xy_s)}{\pm 1 + x^2 +
y^2}ds'\right]\exp\left[-
i\arctan\left(\frac{y_s(0)}{x_s(0)}\right)\right]\right\}_{s = \pm
L},
\end{align*}

\noindent resulting
\begin{align}
z_s(\pm L, t) = \left\{\frac{1 \pm |z|^2}{2}\frac{c_0}{\sqrt
t}e^{i\frac{s^2}{4t}} \exp\left[-i\int_0^s\frac{2(yx_s - xy_s)}{\pm
1 + x^2 + y^2}ds'\right]\exp\left[i\arctan\left(\frac{y_s(0)}{
x_s(0)}\right)\right]\right\}_{s = \pm L}.
\end{align}

\noindent We can give the boundary condition as
\begin{align*}
z_s^{n + 1}(\pm L) = \left\{\frac{1 \pm |z|^2}{2}\frac{c_0}{\sqrt
t}e^{i\frac{s^2}{4t}} \exp\left[-i\int_0^s\frac{2(yx_s - xy_s)}{\pm
1 + x^2 + y^2}ds'\right]\exp\left[i\arctan\left(\frac{y_s(0)}{
x_s(0)}\right)\right]\right\}_{\substack{s = \pm L \\ t = t^n}},
\end{align*}

\noindent which makes the scheme be first order. To obtain a second
order scheme, we do
\begin{align}
\nonumber z_s^{n + 1}(\pm L) & = 2\left\{\frac{1 \pm
|z|^2}{2}\frac{c_0}{ \sqrt t}e^{i\frac{s^2}{4t}}
\exp\left[-i\int_0^s\frac{2(yx_s - xy_s)}{\pm 1 + x^2 +
y^2}ds'\right]\exp\left[i\arctan\left(\frac{y_s(0)}{x_s(0)}\right)\right]\right\}_{\substack{s
= \pm L \\ t = t^n}}
\\
& -\left\{\frac{1 \pm |z|^2}{2}\frac{c_0}{\sqrt
t}e^{i\frac{s^2}{4t}} \exp\left[-i\int_0^s\frac{2(yx_s - xy_s)}{\pm
1 + x^2 + y^2}ds'\right]\exp\left[i\arctan\left(\frac{y_s(0)}{
x_s(0)}\right)\right]\right\}_{\substack{s = \pm L \\ t = t^{n -
1}}}.
\end{align}

\noindent Notice that for unperturbed solutions or for solutions
perturbed with an even perturbation, we have
$$
\exp\left[i\arctan\left(\frac{y_s(0)}{x_s(0)}\right)\right] = 1.
$$
We have found the scheme to be unconditionally stable, and that the
results in the backward case are as good as those obtained with the
boundary condition described in section \ref{sec:self-similar-bc}.
What is more important, this boundary condition seems to work
extraordinarily well for the progressive case. In Figure \ref{figure
frontera sofisticada forward}, considering $t^0 = 1$, we have solved
the equation up to $t = 2$. The results in this case are excellent.
\begin{figure}[!ht]
\centering
\resizebox{9cm}{5cm}{\includegraphics{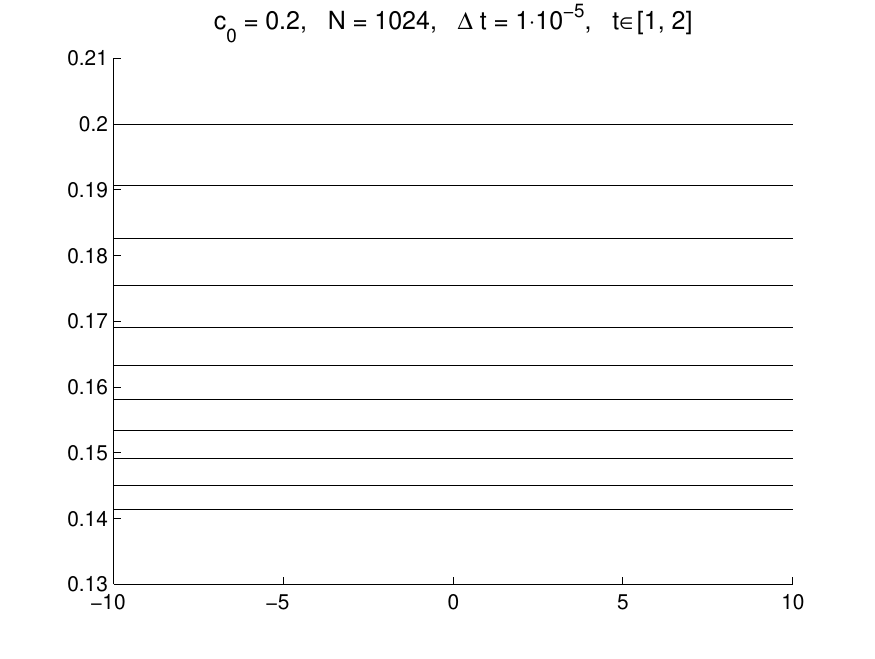}}
{\caption{\small Forward case, using the third boundary conditions.
In comparison with figure \ref{figure diffin forward}, the results
are excellent.} \label{figure frontera sofisticada forward}}
\end{figure}

\subsection{Forward case in time}

The pseudo-spectral method, together with the last boundary
condition, allows to give a complete treatment of the forward case,
unlike with the finite differences.

We consider the following initial value problem:
\begin{align*}
\begin{cases}
z_t = iz_{ss} \mp \dfrac{2i\bar z}{1 \pm |z|^2}z_s^2,
\\
z(s, t) = a^1s\chi_{[0,+\infty)}(s) + a^2s\chi_{(-\infty,0]}(s).
\end{cases}
\end{align*}

\noindent Once  $c_0$ is fixed, we can consider take $a^1$ and $a^2$
as the respective projections of $\mathbf A^1(c_0)$ and $\mathbf
A^2(c_0)$ used in (\ref{2A+- en funcion de c_0}), i.e.,
\begin{align}
a^1 & = \dfrac{\sqrt{\pm(1 - e^{\mp c_0^2\pi})}}{1 + e^{\mp
\frac{c_0^2}{2}\pi}}, & a^2 & = -\dfrac{\sqrt{\pm(1 - e^{\mp
c_0^2\pi})}}{1 + e^{\mp \frac{c_0^2}{2}\pi}}.
\end{align}

\begin{figure}[!ht]
\centering
\resizebox{9cm}{5cm}{\includegraphics{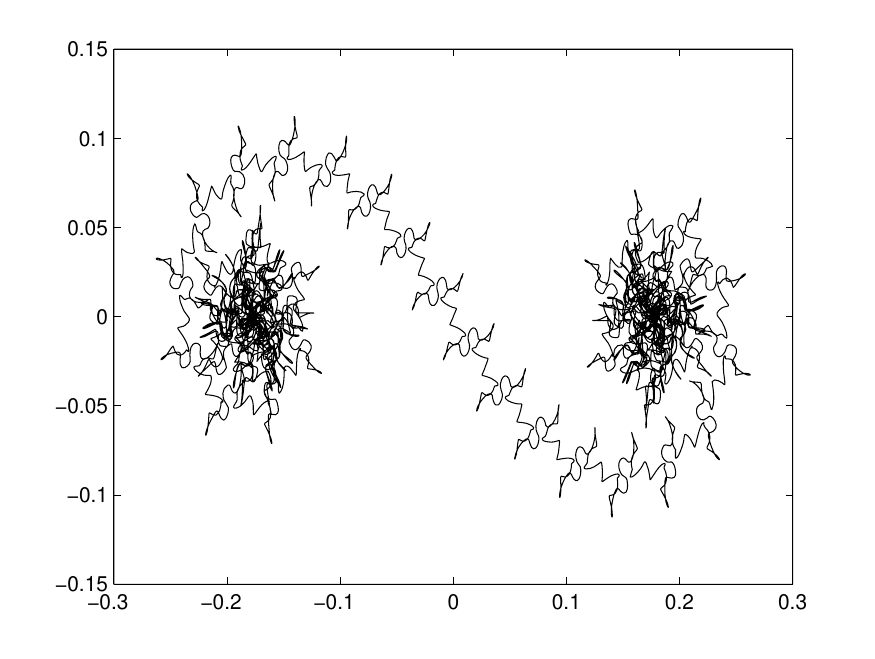}}%
\caption{Fractal creation } \label{figure cheb frac}
\end{figure}

\noindent Therefore, discretizing $[-L, L]$ in $s_i = L\cos(i\pi /
N)$, with $i = 0, \cdots, N$, the numerical initial datum is
\begin{align*}
\begin{cases}
z_i^0 \equiv z(s_i, 0) \equiv a^+, & i\in\{0, \cdots, \frac{N}{2} -
1 \},
\\
z_{N/2}^0 \equiv z(0, 0) \equiv 0,
\\
z_i^0 \equiv z(s_i, 0) \equiv a^-, & i\in\{\frac{N}{2} + 1, \cdots,
N \}.
\end{cases}
\end{align*}
This boundary condition is not adequate for large $t$; when the information reaches the boundary,
it is reflected back into the domain, creating a fractal phenomenon, as illustrated in figure
\ref{figure cheb frac}, where we show the case $c_0 = 0.2$, $N = 16384$, $L = 50$
and $t = 10$.
If, however, we consider for the same parameters a much smaller
time, for instance $t = 0.3$, we observe in Figure \ref{figure cheb
forward1050} that the information, even if it has already reached
the boundary, has only partially rebounded. Thus, in spite of the
great noise at the extremes, if we make a zoom of the central
subinterval $[-10, 10]$, the achieved accuracy is notorious.
\begin{figure}[!ht]
\centering
\resizebox{10cm}{5cm}{\includegraphics{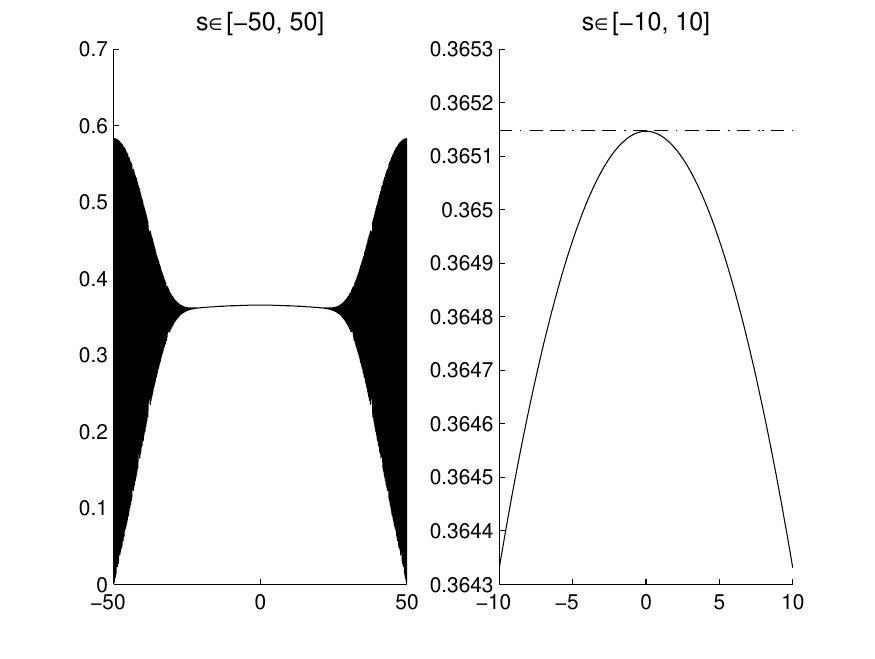}}
{\caption{\small Starting from $t = 0$, with $\Delta t = 10^{-5}$,
at $t = 0.3$ the information has reached the boundary and partially
rebounded. Nevertheless, in $s\in[-10, 10]$, the error of the
curvature is smaller that $9\cdot10^{-4}$, so that portion of curve
is well suited to be our new initial datum. } \label{figure cheb
forward1050}}
\end{figure}
Therefore, we will consider the portion of the curve corresponding
to that subinterval as our new initial datum.

Now, we know $z(s, 0.3)$ at $s_i = 50\cos(i\pi/ 16384)$, $i = 0,
\cdots, 16384$, and have to interpolate spectrally $z(s, 0.3)$ in
the new initial nodes $\tilde s_i$, which we have chosen as
$$
\tilde s_i = 10\cos\left(\frac{i\pi}{1024}\right), \qquad i = 0,
\cdots, 1024.
$$

\noindent Therefore, we have to interpolate $z$ in $1024 + 1$ points
belonging to $[-10, 10]$; this can be done directly from the
expression for $z$ as a function of the coefficients $a_k$:
$$
z(\tilde s_i) = \sum_{k = 0}^Na_kT_k\left(\dfrac{\tilde
s_i}{L}\right) = \sum_{k = 0}^Na_k\cos\left(
k\arccos\left(\dfrac{\tilde s_i}{L}\right)\right).
$$
\begin{figure}[!ht]
\centering
\resizebox{10cm}{5cm}{\includegraphics{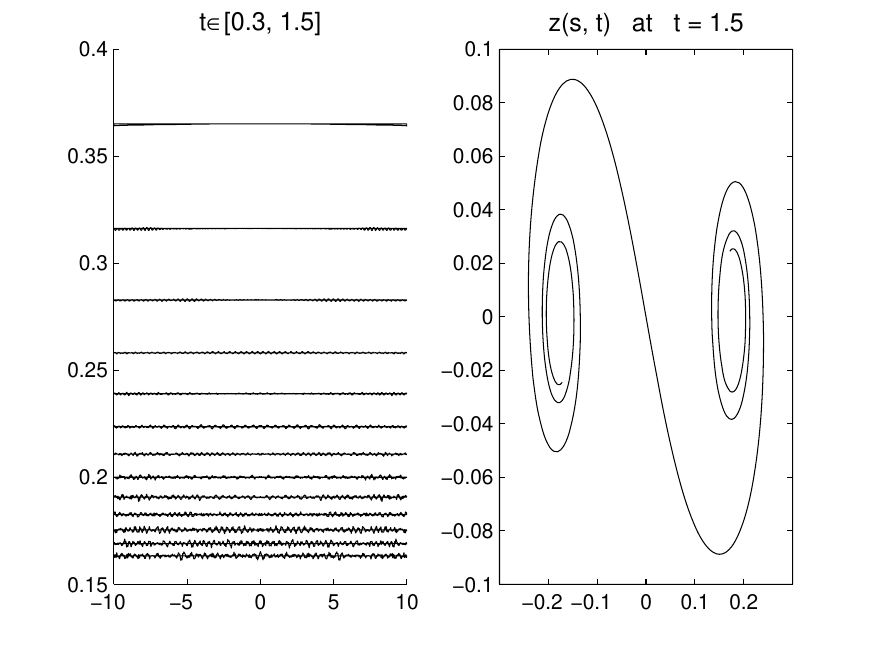}}
{\caption{\small Taking as initial datum the right side of figure
\ref{figure cheb forward1050}, we have advanced till $t = 1.5$, with
$\Delta t = 10^{-5}$ and the third boundary condition. At the right
side, we have drawn the curve obtained at $t = 1.5$.} \label{figur
cheb forwardaprox}}
\end{figure}

\noindent In Figure \ref{figur cheb forwardaprox}, we have
calculated the evolution of $z$ with that new initial datum and the
third boundary condition in $s = \pm 10$. Although there appears
some noise, because the new initial datum is not exact, the results
are acceptable.

To finish, let us mention also that, unlike in the finite difference
case, there are no problems when $\Delta t\to0$ in any of the two
stages.
\section{A stability result}

\label{section stability}

In the previous sections, we have described several schemes to
approximate numerically the self-similar solutions of
\begin{align}
\label{xtposneg2} \Xt & = \Xs\propn\Xss \\
\intertext{and} \label{ttposneg2} \Tt & = \T\propn\Tss,
\end{align}

\noindent which are characterized by
\begin{align}
\label{ctau} c(s, t) = \dfrac{c_0}{2}, \qquad \tau(s, t) =
\dfrac{s}{2t}.
\end{align}

\noindent Since we could not work with all $\mathbb R$, we have
bounded ourselves to $s\in[-L, L]$, considering different boundary
conditions at $s = \pm L$. It is clear that these boundary
conditions, from a numerical point of view, will never be completely
exact; even in the best of the cases there will be some errors,
although they may be as small as the machine precision. Thus, we can
consider all our experiments to be perturbation of the exact
solutions for $\X$, $\T$ and the projection of $\T$, $z$.

In the previous sections, we have seen that the stability of all the
experiments is good, no matter wether we consider very rough
boundary conditions, as fixing $\T(s, t)$ at $\pm L$ or more
elaborated ones. Thus, before concluding this paper, it is
interesting to mention for completeness some recent results by
Banica and Vega \cite{banica1,banica2}, that guarantee theoretically
this stability, at least for small perturbations.

Equations (\ref{xtposneg2}) and (\ref{ttposneg2}) are closely
related to the cubic Schr\"odinger equation or non-linear
Schr\"o\-din\-ger equation (NLS). Indeed, by means of Hasimoto's
transforms
\begin{equation}
\label{hasi} \psi(s, t) = c(s,
t)\exp\left(i\int_0^s\tau(s',t)ds'\right),
\end{equation}

\noindent we obtain the NLS \cite{hasimoto} \cite{buttke2},
\begin{equation}
\label{NLS} i\psi_t + \psi_{ss} \pm \frac{1}{2}[|\psi|^2 + A(t)]\psi
= 0,
\end{equation}

\noindent The term $A(t)$ can be immediately absorbed by means of a
change of variable $\Psi = \psi \exp\left(\mp i/2\int_0^t
A(t')dt'\right)$
\begin{equation}
i\psi_t + \psi_{ss} \pm \frac{1}{2}|\psi|^2\psi = 0.
\end{equation}

\noindent The case with the $+$ sign is known as the focussing case,
the cubic Schr\"odinger equation being denote $\mbox{NLS}^+$; this
corresponds to the Euclidean case of \eqref{xtposneg} and
\eqref{ttposneg}. With the $-$ sign, corresponding to the hyperbolic
case of \eqref{xtposneg} and \eqref{ttposneg}, we have the
defocussing case, the Schr\"odinger cubic being denoted as
$\mbox{NLS}^-$. The Schr\"odinger cubic equation appears in many
contexts \cite{drazin}, as certain non-linear optics phenomena, wave
packets in water and plasma, etc...

Coming back to \eqref{NLS}, in our current problem, from
\eqref{ctau},
\begin{align}
\label{psi(s,t)} \psi(s, t) = \frac{c_0}{\sqrt t}e^{is^2/4t}.
\end{align}

\noindent If we choose $A(t) = -\frac{c_0^2}{t}$, we get a solution
of (\ref{NLS}), with $\psi(0,s) = \sqrt ic_0\delta$, being $\delta$
Dirac's distribution. Thus, going from $\X$ and $z$ to $\psi$ means
to trivialize the problem in a certain way, since, for every $t >
0$, we know the explicit equation.

Now, the solutions (\ref{psi(s,t)}) of the Schr\"odinger cubic
equation (\ref{NLS}) that we are considering satisfy trivially
$$
\int_{\mathbb{R}}|\psi(s, t)|^2ds = \int_{\mathbb{R}}|\psi(0,
s)|^2ds.
$$

\noindent Therefore, from this point of view, our solutions have
infinite energy. Nevertheless, Banica and Vega \cite{banica1,banica2}
have proved that, under certain renormalizations, the solutions have
finite energy. Indeed, starting from
$$
i\psi_t(s, t) + \psi_{ss}(s, t) - \dfrac{1}{2}\left[|\psi(s, t)|^2 -
\dfrac{c_0^2}{t} \right]\psi(s, t) = 0,
$$

\noindent we make the change $u(s, t) = \psi(\sqrt 2s, 2t)$ in order
to absorb the constant $\frac{1}{2}$, obtaining
\begin{align}
iu_t(s, t) + u_{ss}(s, t)\pm\left(|u(s, t)|^2 -
\dfrac{c_0^2}{2t}\right)u(s, t) = 0.
\end{align}

\noindent Applying to this last expression the following conformal
transform
\begin{align} u(s, t) = Tv(s, t) =
\dfrac{e^{i\frac{s^2}{4t}}}{t^{1/2}} v\left(\frac{s}{t},
\frac{1}{t}\right),
\end{align}

\noindent and evaluating in $(s, t) = \left(\frac{s}{t},
\frac{1}{t}\right)$, we get
\begin{align}
v_t(s,t) = -iv_{ss}(s,t) \mp \dfrac{i}{t}\left(\left|v(s,t)
\right|^2 - \frac{c_0^2}{2}\right) v(s,t).
\end{align}

\noindent Therefore, when $0 < t < t_0$, $u$ is solution of
\begin{equation}
\label{0u_t} \left\{
\begin{array}{rcl}
iu_t(s, t) + u_{ss}(s, t)\pm\left(|u(s, t)|^2 -
\dfrac{c_0^2}{2t}\right)u(s, t) & = & 0, \\
u(t_0, s) & = & \dfrac{c_0}{\sqrt{2t_0}}e^{i\frac{s^2}{4}t_0} +
u_1(s),
\end{array}\right.
\end{equation}

\noindent if and only if $v$ is a solution, when $1 / t_0 < t <
\infty$, of
\begin{align}
\label{0v_t} \left\{
\begin{array}{rcl}
v_t(s,t) & = & -iv_{ss}(s,t) \mp \dfrac{i}{t}\left(\left|v(s,t)
\right|^2 - \dfrac{c_0^2}{2}\right) v(s,t),
\\
v(s, 1 / t_0) & = & \dfrac{c_0}{2} + v_0,
\end{array}\right.
\end{align}

\noindent with $v_0(s) = T^{-1}u_1(s)$. Notice that, due to the
conformal transform, the case $t\to0^+$ gets transformed into
$t\to+\infty$. Moreover, our exact unperturbed solutions are now
constant solutions of $u$.

There is an energy naturally associated to (\ref{0v_t}),
\begin{align}
\label{EE} E(t) = \dfrac{1}{2}\int|v_s(s, t)|^2ds \mp
\dfrac{1}{4t}\int\left(|v(s, t)|^2-\dfrac{c_0^2}{2}\right)^2ds.
\end{align}

\noindent Therefore, if $v$ is a solution of (\ref{0v_t}), we get
\begin{align}
\label{0E_t} \dfrac{\partial}{\partial t}E(t) \mp
\dfrac{1}{4t^2}\int\left(|v(s, t)|^2-\dfrac{c_0^2}{2}\right)^2ds =
0.
\end{align}

\noindent This last equation implies, in the defocussing situation
corresponding to the hyperbolic case, that the energy does not grow
when $t\to\infty$. As a consequence of (\ref{0E_t}), Banica and Vega
proved in \cite{banica1} the following theorem for the defocussing
case:
\begin{theorem}

\label{teorema valeria}

For every $t_0 > 0$ and for every $v_0\in\mathcal H^1$, there exists
a unique solution of the initial value problem (\ref{0v_t}), with
$$
v - \dfrac{c_0}{2}\in\mathcal C((1/t_0, \infty), \mathcal H^1).
$$

\end{theorem}

\noindent Banica and Vega also proved that
$$
\liminf_{t\to\infty}\dfrac{1}{t}\int\left(|v(s, t)|^2 -
\dfrac{c_0^2}{2}\right)^2ds = 0,
$$

\noindent which implies that $u$, in the defocussing case, when $0 <
t < t_0$, satisfies
$$
\liminf_{t\to0}\|t|u(t)|^2 - c_0^2\|_2 = 0.
$$

\noindent This last expression is already a stability result for the
singular solution $\dfrac{c_0}{\sqrt{2t_0}}e^{i\frac{s^2}{4}t_0}$ of
(\ref{0u_t}).

Finally, in the defocussing case, writing in terms of the geometric
quantities $c$ and $\tau$ the corresponding energy to \eqref{hasi},
which is a solution of (\ref{NLS}), we get
\begin{align}
\nonumber\noindent\widetilde E(t) & = \frac{t^2}{
4\sqrt2}\int_{-\infty}^{+\infty}\bigg(c_s^2(s, t) + c^2(s,
t)\left(\frac{s}{2t} - \tau(s, t)\right)^2\bigg)
ds \\
& + \frac{1}{16\sqrt 2}\int_{-\infty}^{+\infty}[t c^2(s, t) -
c_0^2]^2ds,
\end{align}

\noindent and, hence,
\begin{equation}
\frac{d}{dt}\widetilde E(t) - \frac{1}{16\sqrt 2\ t}
\int_{-\infty}^{+\infty} [t c^2(s, t) - c_0^2]^2ds = 0
\end{equation}

\noindent and
\begin{equation}
\liminf_{t\to0}\|t|c|^2-c_0^2\|_2 = 0.
\end{equation}

\noindent Observe that our solutions are precisely such that
$\widetilde E(t) = 0$, for all $t > 0$.

In the focussing setting (i.e. the case of the sphere) the stability is
much more delicate. In \cite{banica2}, and under a smallness assumption
in the curvature $c_0$, they construct  a global solution such that
$E(t)$ given in \eqref{EE} is finite. The main difficulty comes from
the long range character of the non-linear potential that appears in
\eqref{0v_t}. This  implies the existence of a logarithmic phase that
has as a consequence the non existence of the limit at infinity for the
solutions. However they also prove, this time in the focussing case,
that this logarithmic divergence disappears when the tangent vector $T$
is computed, so that there is stability for $T$. Remember that in order
to compute $T$ one has to integrate once the curvature and twice the
torsion. Therefore thanks to the oscillations the integrals converge
without any difficulty.

\section{Conclusions}

In this paper, we have tried to reproduce numerically the behavior
of the self-similar solutions of
\begin{align}
\label{xtposneg3} \Xt & = \Xs\propn\Xss \\
\intertext{and} \label{ttposneg3} \Tt & = \T\propn\Tss,
\end{align}

\noindent which develop a singularity at finite time. These
solutions, characterized by
\begin{align}
c(s, t) = \dfrac{c_0}{2}, \qquad \tau(s, t) = \dfrac{s}{2t},
\end{align}

\noindent form a one-parameter family, where $c_0$ is precisely the
family parameter.

The singularities happen at $t = 0$, going backwards in time, but,
since both \eqref{xtposneg3} and \eqref{ttposneg3} are
time-reversible, we could consider an equivalent problem, where we
advance in time and the singularity happens.

In Section \ref{section findif}, we have given a finite-difference
scheme to study the self-similar solutions of $\T$. $\T$ determines
completely $\X$, except for a constant, determined by $\X(0, t) =
2c_0\sqrt t(0, 1, 0)$. We have studied mainly the backward case:
starting from $t = 1$, we have tried to reproduce the formation of
the singularity. Since we cannot consider all $\mathbb R$, we have
bounded ourselves to $s\in[-L, L]$, being necessary to give boundary
conditions at $s=\pm L$. From the asymptotics of Theorem
\ref{teorema} for $\T(s, 1)$, we have deduced two approximated
boundary conditions, depending on whether we choose the leading term
of  $\T(s, 1)$, or the first two terms in the expansion.

Considering $\T(s, 1)$ constant at $s = \pm L$ gives us good results
from a qualitative point of view: the energy is preserved with
several precision digits and we get also the approximated value of
$c(0, t)$, even for small times. Thus, the finite energy tends to
concentrate on the origin, approximating the formation of the
singularity at $t = 0$. In the exact solution, the energy was
infinite for all $t$, but the behavior is the same: all the energy
tends to concentrate at $s = 0$.

Since the bigger is $L$, the smaller are the times for which we
recover $c(0, t)$, there seems to be evidence that we could
approximate the exact problem by making $L$ tend to infinity. It
would be very interesting to prove this analytically.

The second boundary condition is obtained considering the first
non-constant term in the asymptotics of $\T(s, 1)$, and that $\T(s,
t) = \T(s / \sqrt t, 1)$. For not too big $L$ and provided that
$\Delta s$ is small enough, it allows us to recover the solutions
with big accuracy even for small times; moreover, not only for $s =
0$, but for all $s$. Nevertheless, since $|\Delta t| = \mathcal
O(\Delta s^2)$, this is quite expensive from a computational point
of view.

When $L$ is bigger, a uniform distribution of the nodes is not
adequate, because there is a lack of resolution for $s$ near the
boundary. That suggests that we use an alternative node
distribution, given by the Chebyshev nodes, which are distributed in
a much more suitable way. Since, with an explicit scheme, we would
have now a $|\Delta t| = \mathcal O(N^{-4})$ restriction, we have
projected stereographically $\T$ over $\mathbb C$, obtaining
\begin{align}
\label{z_t4} z_t = iz_{ss} \mp \frac{2i\bar z}{1 \pm |z|^2}z_s^2.
\end{align}

\noindent We have implemented an implicit-explicit method in time,
with a pseudo-spectral method in space, considering different
boundary conditions. This method has several interesting advantages:
its stability is much better for a given $N$; it allows considering
much bigger $L$; it allows implementing boundary conditions that do
not need external information and $\X$ can be recovered immediately
with spectral accuracy. Since designing an adaptive version, i.e.,
refining the grid when necessary, is natural and straight-forward,
the method is globally much more efficient. Lastly, it allows to
make a full treatment of the forward case: starting from a singular
datum at $t = 0$, we can recover the solutions in two stages, by
adapting at $t = \varepsilon$ adequate boundary conditions.

Finally, let us mention that all the experiments we have done give
evidence of the stability of equations \eqref{xtposneg3},
\eqref{ttposneg3} and \eqref{z_t4} from a numerical point of view.
Indeed, when giving the boundary conditions at $s = \pm L$, it is
impossible to do it exactly and we will always be introducing some
perturbations; when, for instance, we fix $\T(s, t)$ to be constant
at $s = \pm L$, these perturbations will be quite big. Even if we
calculated the exact value of $\X$, $\T$ or $z$ at $s = \pm L$,
there would always be a tiny error attributable to machine
precision. Therefore, in Section \ref{section stability}, we have
mention a recent stability result by Banica and Vega, concerning the
self-similar solutions we have studied. This offers, somehow, a
theoretical support to this paper.

\section{Acknowledgements}

The authors would want to express their gratitude to M.~A. Fontelos
and E. Zuazua, for allowing them to use the Odisea cluster, located
at the Department of Mathematics of the Universidad Aut\'onoma de
Madrid, to perform extensive numerical simulations. Part of this
work was carried out while CJGC was visiting the Univerdad del
Pa{\'\i}s Vasco in Bilbao. CJGC is grateful for the invitation, and
the hospitatility of all the faculty and staff at the university.

\bibliographystyle{siam}


\end{document}